\documentclass[11pt]{amsart}
\usepackage{latexsym,graphicx}
\numberwithin{equation}{section}
\theoremstyle{plain}

\theoremstyle{remark}

\theoremstyle{definition}

\newcommand{\D}{{\mathcal D}}
\newcommand{\E}{\mathcal E}

\newcommand{\G}{{\mathcal G}}

\renewcommand{\H}{\mathbb H}

\newcommand{\K}{{\mathcal K}}
\renewcommand{\L}{{\mathcal L}}
\newcommand{\M}{{\mathcal M}}
\newcommand{\N}{\mathbb N}

\newcommand{\V}{{\mathcal V}}

\newcommand{\dist}{\operatorname{dist}}

\newcommand{\fp}{\operatorname{FP}}

\newcommand{\hs}{\operatorname{HS}}

\newcommand{\Int}{\operatorname{Int}}

\renewcommand{\span}{\operatorname{span}}

\newcommand{\supp}{\operatorname{Supp}}

\def\half{{1 \over 2}}

\newcommand{\oa}{\overrightarrow}

\newcommand{\ol}{\overline}

\def\XXint#1#2#3{{\setbox0=\hbox{$#1{#2#3}{\int}$}
      \vcenter{\hbox{$#2#3$}}\kern-.5\wd0}}

\def\elll{\ell}

\begin{document}

\def\cal{\mathcal}

\font\tpt=cmr10 at 12 pt
\font\fpt=cmr10 at 14 pt

\font \fr = eufm10

%\font\AAA=Times.dfont  at 12pt
 %\font\BBB=Times.dfont at 8pt

%\font\AAA=cmr10 at 12pt
%\font\BBB=cmr10 at 8pt

%\def\AAA{\bf}
%\def\BBB{\bf}

\overfullrule=0in

\def\boxit#1{\hbox{\vrule
 \vtop{%
  \vbox{\hrule\kern 2pt %
     \hbox{\kern 2pt #1\kern 2pt}}%
   \kern 2pt \hrule }%
  \vrule}}

  \def\harr#1#2{\ \smash{\mathop{\hbox to .3in{\rightarrowfill}}\limits^{\scriptstyle#1}_{\scriptstyle#2}}\ }

\def\AA{1}
\def\BB{2}
\def\CC{3}
\def\DD{4}
\def\EE{5}
\def\FF{6}
\def\GGG{7}
\def\HH{8}
\def\II{9}
\def\JJ{10}
\def\KK{11}
\def\LL{12}
\def\MM{13}

\def\ALL{1}
\def\BTA{2}
\def\BL{3}
\def\BRE{4}
\def\CNS{5}
\def\CIL{6}
\def\CRA{7}
\def\DDD{8}
\def\DDR{9}
\def\GEO{10}
\def\HYP{11}
\def\BEL{12}
\def\AC{13}
\def\SURVEY{14}
\def\NOTES{15}
\def\AET{16}
\def\LAG{17}
\def\KRY{18}
\def\PLI{19}
\def\RT{20}
\def\SLO{21}
\def\TRU{22}
\def\TWC{23}
\def\WAL{24}

 \def\GG{{{\bf G} \!\!\!\! {\rm l}}\ }

\def\GL{{\rm GL}}

\def\bll{I \!\! L}

\def\IFF{\qquad\iff\qquad}
\def\bra#1#2{\langle #1, #2\rangle}
\def\bbf{{\bf F}}
\def\bbj{{\bf J}}
\def\Jtn{{\bbj}^2_n}  \def\JtN{{\bbj}^2_N}  \def\JoN{{\bbj}^1_N}
\def\jt{j^2}
\def\jtx{\jt_x}
\def\Jt{J^2}
\def\Jtx{\Jt_x}
\def\bpp{{\bf P}^+}
\def\bpt{{\wt{\bf P}}}
\def\fsh{$F$-subharmonic }
\def\mo{monotonicity }
\def\jet{(r,p,A)}
\def\ss{\subset}
\def\sse{\subseteq}
\def\half{\hbox{${1\over 2}$}}
\def\smfrac#1#2{\hbox{${#1\over #2}$}}
\def\oa#1{\overrightarrow #1}
\def\dim{{\rm dim}}
\def\dist{{\rm dist}}
\def\codim{{\rm codim}}
\def\deg{{\rm deg}}
\def\rank{{\rm rank}}
\def\log{{\rm log}}
\def\Hess{{\rm Hess}}
\def\Hessyp{{\rm Hess}_{\rm SYP}}
\def\trace{{\rm trace}}
\def\tr{{\rm tr}}
\def\max{{\rm max}}
\def\min{{\rm min}}
\def\span{{\rm span\,}}
\def\Hom{{\rm Hom\,}}
\def\det{{\rm det}}
\def\End{{\rm End}}
\def\Sym{{\rm Sym}^2}
\def\diag{{\rm diag}}
\def\pt{{\rm pt}}
\def\Spec{{\rm Spec}}
\def\pr{{\rm pr}}
\def\Id{{\rm Id}}
\def\Grass{{\rm Grass}}
\def\Herm#1{{\rm Herm}_{#1}(V)}
\def\arr{\longrightarrow}
\def\supp{{\rm supp}}
\def\Link{{\rm Link}}
\def\Wind{{\rm Wind}}
\def\Div{{\rm Div}}
\def\vol{{\rm vol}}
\def\foral{\qquad {\rm for\ all\ \ }}
\def\fpsh{{\cal PSH}(X,\f)}
\def\Core{{\rm Core}}
\def\dis{f_M}
\def\Re{{\rm Re}}
\def\rn{\bbr^n}
\def\pp{\cp^+}
\def\plp{\cp_+}
\def\Int{{\rm Int}}
\def\cix{C^{\infty}(X)}
\def\Gr#1{G(#1,\rn)}
\def\Symn{{\Sym(\rn)}}
\def\SymN{{\Sym(\bbr^N)}}
\def\Gpn{G(p,\rn)}
\def\fd{{\rm free-dim}}
\def\SA{{\rm SA}}
 \def\cd{{\cal C}}
 \def\cdt{{\widetilde \cd}}
 \def\cm{{\cal M}}
 \def\cmt{{\widetilde \cm}}

\def\Theorem#1{\medskip\noindent {\bf THEOREM \bf #1.}}
\def\Prop#1{\medskip\noindent {\bf Proposition #1.}}
\def\Cor#1{\medskip\noindent {\bf Corollary #1.}}
\def\Lemma#1{\medskip\noindent {\bf Lemma #1.}}
\def\Remark#1{\medskip\noindent {\bf Remark #1.}}
\def\Note#1{\medskip\noindent {\bf Note #1.}}
\def\Def#1{\medskip\noindent {\bf Definition #1.}}
\def\Claim#1{\medskip\noindent {\bf Claim #1.}}
\def\Conj#1{\medskip\noindent {\bf Conjecture \bf    #1.}}
\def\Ex#1{\medskip\noindent {\bf Example \bf    #1.}}
\def\Qu#1{\medskip\noindent {\bf Question \bf    #1.}}
\def\Exercise#1{\medskip\noindent {\bf Exercise \bf    #1.}}

\def\HoQu#1{ {\AAA T\BBB HE\ \AAA H\BBB ODGE\ \AAA Q\BBB UESTION \bf    #1.}}

\def\pf{\medskip\noindent {\bf Proof.}\ }
\def\qed{\hfill  $\vrule width5pt height5pt depth0pt$}
\def\equdef{\buildrel {\rm def} \over  =}
\def\qedqed{\hfill  $\vrule width5pt height5pt depth0pt$ $\vrule width5pt height5pt depth0pt$}
\def\mathqed{  \vrule width5pt height5pt depth0pt}

\def\V{W}

\def\df{d^{\phi}}
\def\hk{\_{\rm l}\,}
\def\n{\nabla}
\def\w{\wedge}

\def\cu{{\cal U}}   \def\cc{{\cal C}}   \def\cb{{\cal B}}  \def\cz{{\cal Z}}
\def\cv{{\cal V}}   \def\cp{{\cal P}}   \def\ca{{\cal A}}
\def\cw{{\cal W}}   \def\co{{\cal O}}
\def\ce{{\cal E}}   \def\ck{{\cal K}}
\def\ch{{\cal H}}   \def\cm{{\cal M}}
\def\cs{{\cal S}}   \def\cn{{\cal N}}
\def\cd{{\cal D}}
\def\cl{{\cal L}}
\def\cp{{\cal P}}
\def\cf{{\cal F}}
\def\ccr{{\cal  R}}

\def\gerG{{\fr{\hbox{g}}}}
\def\gerB{{\fr{\hbox{B}}}}
\def\gerR{{\fr{\hbox{R}}}}
\def\p#1{{\bf P}^{#1}}
\def\vf{\varphi}

\def\wt{\widetilde}
\def\wh{\widehat}

\def\and{\qquad {\rm and} \qquad}
\def\arr{\longrightarrow}
\def\ol{\overline}
\def\bbr{{\mathbb R}}\def\bbh{{\mathbb H}}\def\bbo{{\mathbb O}}
\def\bbc{{\mathbb C}}
\def\bbq{{\mathbb Q}}
\def\bbz{{\mathbb Z}}
\def\bbp{{\mathbb P}}
\def\bbd{{\mathbb D}}

\def\a{\alpha}
\def\b{\beta}
\def\d{\delta}
\def\e{\epsilon}
\def\f{\phi}
\def\g{\gamma}
\def\k{\kappa}
\def\l{\lambda}
\def\o{\omega}

\def\s{\sigma}
\def\x{\xi}
\def\z{\zeta}

\def\D{\Delta}
\def\L{\Lambda}
\def\G{\Gamma}
\def\O{\Omega}

\def\bd{\partial}
\def\bdf{\partial_{\f}}
\def\lag{Lagrangian}
\def\psh{plurisubharmonic }
\def\ph{pluriharmonic }
\def\pph{partially pluriharmonic }
\def\omp{$\omega$-plurisubharmonic \ }
\def\ffl{$\f$-flat}
\def\PH#1{\widehat {#1}}
\def\lloc{L^1_{\rm loc}}
\def\dbar{\ol{\partial}}
\def\lp{\Lambda_+(\f)}
\def\lpp{\Lambda^+(\f)}
\def\bo{\partial \Omega}
\def\Ob{\overline{\O}}
\def\fc{$\phi$-convex }
\def\PSH{{ \rm PSH}}
\def\SH{{\rm SH}}
\def\totr{ $\phi$-free }
\def\BM{\lambda}
\def\Der{D}
\def\CH{{\cal H}}
\def\RH{\overline{\ch}^\f }
\def\pconv{$p$-convex}
\def\MA{MA}
\def\lagpsh{Lagrangian plurisubharmonic}
\def\hermsk{{\rm Herm}_{\rm skew}}
\def\PSHl{\PSH_{\rm Lag}}
 \def\ppsh{$\pp$-plurisubharmonic}
\def\fp{$\pp$-plurisubharmonic }
\def\fh{$\pp$-pluriharmonic }
\def\Symn{\Sym(\rn)}
 \def\ci{C^{\infty}}
\def\USC{{\rm USC}}
\def\LSC{{\rm LSC}}
\def\fa{{\rm\ \  for\ all\ }}
\def\ppc{$\pp$-convex}
\def\cpt{\wt{\cp}}
\def\ft{\wt F}
\def\ob{\overline{\O}}
\def\Be{B_\e}
\def\K{{\rm K}}

\def\M{{\bf M}}
\def\N#1{C_{#1}}
\def\ds{Dirichlet set }
\def\dir{Dirichlet }
\def\Fa{{\oa F}}
\def\TR{{\cal T}}
 \def\ISO{{\rm ISO_p}}
 \def\Span{{\rm Span}}

\def\II{1}
\def\AB{2}
\def\AA{3}
\def\BB{4}
\def\CC{5}
\def\DD{6}
\def\EE{7}
\def\FF{8}

\vskip .4in

\def\E{E}
\def\fpsi{{F_f(\psi)}}
\def\bL{{\bf \Lambda}}
\def\abbf{\oa\bbf}
\def\pitwo{{\pi \over 2}}
\def\Fth{{\bbf_\theta}}
\def\aFth{{\abbf_{\theta}}}
\def\AI{{\oa{\Int} \, \bbf}}
\def\AIth{{\oa{\Int} \, \Fth}}
\def\AIthk{{\oa{\Int} \,{\bbf_{\theta_k}}}}
\def\ggg{{\mathfrak g}}

\font\headfont=cmr10 at 14 pt
\font\aufont=cmr10 at 11 pt

\title[THE SPECIAL LAGRANGIAN POTENTIAL EQUATION]
{\headfont PSEUDOCONVEXITY FOR THE SPECIAL LAGRANGIAN POTENTIAL EQUATION}

\date{\today}
\author{ F. Reese Harvey and H. Blaine Lawson, Jr.}
\thanks
{Partially supported by the NSF}

\maketitle

\begin{abstract}
The Special Lagrangian Potential Equation for a function $u$  on a domain $\O\ss\rn$ is given by 
$\tr\{\arctan(D^2 \,u) \} = \theta$  for a contant $\theta \in (-n {\pi\over 2}, n {\pi\over 2})$.
For $C^2$ solutions the graph of $Du$ in $\O\times\rn$ is a special Lagrangian submanfold.
Much has been understood about the Dirichlet problem for this equation, but the existence result
relies on explicitly computing the associated boundary conditions (or, otherwise said,  computing the pseudo-convexity for the associated potential theory).  This is done in this paper, and the answer is interesting.
The result carries over to many related equations --  for example, those obtained by taking 
$\sum_k \arctan\, \l_k^\ggg = \theta$ where $\ggg : \Symn\to \bbr$ is a G\aa rding-Dirichlet polynomial
which is hyperbolic with respect to the identity.  A particular example of this is the 
deformed Hermitian-Yang-Mills equation which appears in mirror symmetry.
Another example is $\sum_j \arctan \kappa_j = \theta$ where $\kappa_1, ... , \kappa_n$ are the 
principal curvatures of the graph of $u$ in $\O\times \bbr$.

We also discuss the inhomogeneous Dirichlet Problem 

\centerline{ $\tr\{\arctan(D^2_x \,u)\} = \psi(x)$}

\noindent
where $\psi : \ob\to (-n {\pi\over 2}, n {\pi\over 2})$.  This equation has the feature that
the pull-back of $\psi$ to the Lagrangian submanifold  $L\equiv {\rm graph}(Du)$ is the 
phase function  $\theta$ of the tangent spaces of $L$.  On $L$   it  satisfies the equation 
$\nabla \psi = -JH$  where $H$ is the mean curvature vector field of $L$.

\end{abstract}

%{\small\tableofcontents}

\vskip .1in

 \vfill\eject

%\section{Introduction}
%\label{intro}
\ \ 
\vskip 1in

\centerline{\bf Table of Contents} \bigskip

%{{\parindent= .1in\narrower 

 \hskip .5 in  \II. Introduction

 \hskip .5in    \AB.  Geometric Conditions for Strict $\abbf_\theta$-Boundary Convexity

 \hskip .5 in  \AA. Preliminaries 
 
 \hskip 1in The Special Lagrangian Potential Equation $\bbf_\theta$
 
  \hskip 1in   The Asymptotic Interior
  
   \hskip 1in Branches
   
    \hskip 1in Some General Results on the Pure Second-Order 
    \\ \medskip  \hskip 1.5in  Dirichlet Problem

  \hskip .5 in  \BB.   Computing the Asymptotic Interior of $\bbf_\theta$

  \hskip .5 in  \CC.   The Refined  Dirichlet Problem for $\bbf_\theta$ % the SL Potential Equation of 
%    \\ \medskip  \hskip 1.5in  Phase $\theta$.

    \hskip .5 in  \DD.   The Inhomogeneous Dirichlet Problem for the  SL Potential 
    
\hskip .7 in    Operator

    \hskip .5 in  \EE.  A Generalized Version of the Main Theorem \BB.1

    \hskip .5 in  \FF.  Results on Riemannian Manifolds
 
 \vskip .2in 
Appendix A.  A Geometric Interpretation of the Inhomogeneous DP.

Appendix B.  Remarks Concerning Convexity

\vfill\eject

%%%%%%%%%%%%%%%%%%%%%%%%%%%%%%%%%%%%%%%%%%%%%%%%
%%%%%%%%%%%%%%%%%%%%%%%%%%%%%%%%%%%%%%%%%%%%%%%%
%%%%%%%%%%%%%%%%%%%%%%%%%%%%%%%%%%%%%%%%%%%%%%%%
%%%%%%%%%%%%%%%%%%%%%%%%%%%%%%%%%%%%%%%%%%%%%%%%
%%%%%%%%%%%%%%%%%%%%%%%%%%%%%%%%%%%%%%%%%%%%%%%%
%%%%%%%%%%%%%%%%%%%%%%%%%%%%%%%%%%%%%%%%%%%%%%%%
%%%%%%%%%%%%%%%%%%%%%%%%%%%%%%%%%%%%%%%%%%%%%%%%

\noindent {\headfont \II.    Introduction.}

%Existence for this Dirichlet problem requires computing the asymptotic cone for  the subequation $\bbf_\theta$.
%%For $\theta > (n-1){\pi\over 2}$ this was done in [CNS].  The main point of this article  
%is to compute this asymptotic cone for all $\theta$, thereby providing the widest class of domains
%$\O$ where existence holds.

 The special Lagrangian potential equation was introduced in [CG] back in 1982.
 Its solutions $u$ were shown to have the property that the graph $p=\nabla u$ in $\rn\times\rn = \bbc^n$
is a Lagrangian submanifold which is  absolutely volume-minimizing, 
and the linearization at any solution is elliptic.    Many examples of
these {\sl Special Lagrangian submanifolds} were given in [CG], but the Dirichlet problem for this equation
was a difficult challenge.  It was first  solved in the $C^\infty$-case, for  ``large phases'',  in the beautiful paper
of Caffarelli, Nirenberg and Spruck [CNS].    Then for all phases, existence and uniqueness of viscosity solutions   in the $C^0$ case were established  in [DD].

More specifically, this equation, with {\bf phase} $\theta\in (-n{\pi\over 2}, n{\pi\over 2})$, is written:
$$
 f(D^2 u)\ \equdef\           \tr\left\{ \arctan (D^2 u)\right\}\ =\ \theta.
\eqno{(\II.1)}
$$
The associated special Lagrangian submanifolds have the property that  the $n$-form  
$\Phi_\theta \equiv {\rm Im}\{e^{-i\theta} dz_1\wedge \cdots \wedge dz_n\}$ vanishes identically
on them,  and (with appropriate orientation) they are  calibrated by 
${\rm Re}\{e^{-i\theta} dz_1\wedge \cdots \wedge dz_n\}$.
Now it is an important fact that the set of  algebraic solutions, i.e., the constraint set on second
derivatives, 
$$
{\mathcal H}_\theta\ \equiv\  \left \{ A\in\Symn :  \Phi_\theta \bigr|_{\rm graph(A)}=0 \right\}
\eqno{(\II.2)}
$$
is not connected.  Specifically, the equation (\II.1) gives rise to the subequation
$$
\bbf_\theta \equiv \{A\in\Symn : \tr\left\{ \arctan (A)\right\}\geq\theta\},
\eqno{(\II.3)}
$$
and we have that 
$$
{\mathcal H}_\theta \ =\  \bigcup_{\theta' \  \equiv\ \theta \ {\rm mod}  \, 2\pi}  \partial \bbf_{\theta'}
\eqno{(\II.4)}
$$
where each of the equations $\partial \bbf_{\theta'}$ is a connected component.

Now the solutions to  (\II.1) in [CNS] were for phases  $\theta$ with 
$(n-1){\pi\over 2} <  |\theta| < n{\pi\over 2}$.  These are the phases where the operator is concave
for $\theta>0$ and convex for $\theta <0$.
In [DD] the solutions are obtained  for all phases, i.e.,  $\theta\in (-n{\pi\over 2}, n{\pi\over 2})$.

The  best answer for the existence question for this Dirichlet problem requires computing 
the asymptotic cone for  the subequation 
$ \bbf_\theta$.
For $|\theta| > (n-1){\pi\over 2}$ this was done in [CNS].  
The unfinished business, which is completed in this article,  
is to compute this asymptotic cone for all $\theta$, thereby providing the widest class of domains
$\O$ where existence holds, or, said differently, providing the appropriate notion of boundary pseudo-convexity
for the potential theory associated to the SL-operator (\II.1).  
As explained below, the appropriate notion of pseudo-convexity for $\bbf_\theta$ only depends 
on $|\theta|$.

Interestingly, as the phases get closer to zero, the restriction on pseudo-convexity gets weaker and 
weaker.  Therefore, for the various phases $\theta'$ appearing in (\II.4) above, existence of solutions
to the Dirichlet Problem holds on larger and larger categories of domains   as $|\theta'|\to 0$.

Specifically, our result is the following Theorem.   There is associated to  $\bbf_\theta$  its {\sl  asymptotic 
subequation} $\abbf_\theta$, which is a cone with vertex at the origin (see Section \AA\ for the definition). 
Consider a domain $\O\ss\rn$ with smooth boundary $\bo$.
Let ${\rm II}_{x,\bo}$ denote the second fundamental form of $\bo$ at $x$ with respect to the interior unit normal 
$n = n_x$, and let $P_n$  be orthogonal projection onto the line $\bbr n$.
Then we say that $\bo$ is {\bf strictly $\bbf_\theta$-convex} at $x$, if 
$$
{\rm II}_{x, \bo} + tP_{n} \ \in \ \Int \oa \bbf_\theta
\ \ {\rm for  } \ t>>0,
$$
(See also (\AA.9b) for an equivalent definition.)

In approaching the Dirichlet problem, one needs to consider the subequation $\bbf_\theta$
 {\bf and its dual} $\wt \bbf_\theta \equiv \sim\{- \Int \bbf_\theta\}$.  An $\bbf_\theta$-subharmonic function  $u$ is a 
 subsolution, and if $-u$ is an  $\wt \bbf_\theta$-subharmonic, then $u$ is a supersolution.  
 If both conditions hold, then $u$ is a viscosity solution to the equation (\II.1).
 For the SL protential equation this duality is very pretty.  The dual of $\bbf_\theta$ is
 $$
 \wt \bbf_\theta \ =\ \bbf_{-\theta}.
 $$

 Now this duality
 carries over to the boundary conditions necessary for existence.
 At each point the boundary $\bo$ of the domain must be both  strictly $\bbf_\theta$-convex
 and  strictly $\wt \bbf_\theta$-convex.  If  $\theta' <\theta$,  then $\bbf_{\theta'} \supset \bbf_\theta$.
Hence, we see that the boundary condition for the $\bbf_\theta$-Dirichlet problem  is that 
$\bo$ must be  strictly $\bbf_{|\theta|}$-convex at every point.
That is, strict $\bbf_{|\theta|}$ convexity of $\bo$ is exactly the condition necessary to establish existence for
the Dirichlet Problem for the equation $\bbf_\theta$ for all continuous boundary data.  So we want to compute 
$\oa \bbf_\theta$ explicitly.

Consider the operator $$f(A) \equiv \tr \{ \arctan (A)\}$$ defined on $\bbf \equiv \Symn$.  
This operator $f$ has values precisely in the open interval $(-n{\pi  \over 2},  n{\pi  \over 2})$.
Of these,  there are $n-1$ 
$$
{\bf Special\ Phases \ (or\  Values):}   \hskip .3in  \theta_k \ =\ (n-2k) {\pi\over 2}, \quad k=1,..., n-1.  
$$
Removing these $n-1$ special values, the remaining set of values is 
the disjoint union of $n$ open
$$
{\bf  Phase \ Intervals:} \ \ 
I_k \ =\  \left ( (n-2k) {\pi\over 2}, \ (n-2(k-1)) {\pi\over 2} \right), \ k=1,...,n.
$$

Let  $\l_1(A)\leq\l_2(A)\leq \cdots$ denote  the ordered eigenvalues of $A$, and let $\s_k(A)$ be the 
$k^{\rm th}$ elementary symmetric function of the eigenvalues of $A$.

\Theorem {\BB.1}  {\sl
The asymptotic subequation $\abbf_\theta$ of\ \  $\bbf_\theta$, for $\theta \in (-n{\pi  \over 2},  n{\pi  \over 2})$,
    is given as follows.

(1) \ \ If $\theta\in I_k$ ($k=1, ... ,n$), then
$$
\abbf_\theta\ =\ \bL_k \ \equiv\ \{A\in\Symn : \l_k(A)\geq 0\}.
$$

(2) \ \ If $\theta_k$ ($k=1, ... ,n-1$) is a special value, then}
$$
\abbf_{\theta_k}\ =\ \bL_k^{\s_{n-1}}.
$$
The set $\bL_k^{\s_{n-1}}$ depends on the G\aa rding polynomial $\s_{n-1}(A)$ whose eigenvalues {\sl cannot}
be computed in terms of the eigenvalues of $A$. (See the subsection ``branches'' in section \AA.)
However, from Proposition \BB.5  we have the following.

\Prop{\II.1} {\sl
For $k=1, ... , n-1$,  the set  $\bL_k^{\s_{n-1}}$ is a disjoint union $\bL_k^{\s_{n-1}} = \bL_k \cup E_k$
where
$$
E_k \ \equiv \ 
 \left(\bL_{k+1} \sim \bL_{k}\right) \cap\left\{  \s_{j}(A)\,\s_{n}(A) < 0\right\}
$$
and where $j$ is defined by the condition: $\s_j(A) \neq 0$ and $\s_\ell(A)=0$ for $j<\ell\leq n-1$.
}

An application of this result is to the existence question for the Dirichlet problem.
For this one needs to know the interior of $\bL_k^{\s_{n-1}}$.  Part (2) of the following
is not obvious from Proposition \II.1.  Note that  the terms with  $j<n-1$ do not enter.
This is done is Proposition \BB.5.

\Theorem {\BB.6} {\sl
The interior of $\bL_k^{\s_{n-1}}$ is given as follows.

(1)\ \ If $\theta\in I_k$ ($k=1, ... ,n$), then
$$
\Int \abbf_\theta\ =\  \Int \bL_k \ \equiv\ \{A\in\Symn : \l_k(A) > 0\}.
$$

 (2) \ \ If $\theta_k$ ($k=1, ... ,n-1$) is a special value, then
$$
\Int \abbf_{\theta_k}\ =\ \{\Int \bL_{k+1} \cap \bL_k\} \cup E_k^\star
$$
where
$$
E_k^\star \ \equiv \ 
 \left(\Int \bL_{k+1} \sim \bL_{k}\right) \cap\left\{  \s_{n-1}(A)\,\s_{n}(A) < 0\right\}
$$
}

Using Theorem \BB.6 we give a deeper version of our general results on the Dirichlet problem [DD]
for this equation.  Note that as $k$ increases from 1 to $[n/2]$, the sets $\bL_k$ get quite large.
The first $k-1$ eigenvalues of $A$ can be arbitrarily negative.  So for the intervals close to the origin
the geometric constraints on the second fundamental form of the boundary of the
domain are quite loose.   This general existence and uniqueness result for the SL Dirichlet problem is given in Section \CC.

There has also been work on the inhomogeneous Dirichlet problem
$$
\tr\left\{ \arctan (D^2 u)\right\}\ =\ \psi(x)
\eqno{(\II.5)}
$$
where $\psi$ is a continuous function on the closed domain with values in a high phase interval.
Solutions in the $C^\infty$ category were obtained by
Tristan  Collins, Sebastien Picard and Xuan Wu [CPW] where the interval is the critical one:
$( (n-2){\pi\over 2}, n{\pi\over 2})$. (See Theorem \DD.3.) The analogue of this result in the continuous case was
done by  S. Dinew, H.-S. Do and T. D. T\^o  in [DDT].  Recently, Marco Cirant and Kevin Payne have established comparison for this equation when $\psi$ does not take on a special value, in other words, when $\psi$ takes its values in (any) one of the phase intervals $I_k$ (see [CP, \S 6.4] and Theorem \DD.2 below).

In fact using [CP, \S 6.4]  and Theorem \AB.1 we get the following.  Let $\O\ss\ss\rn$ be a domain
with smooth boundary, and at each $x\in\O$ let $\k_1(x)\leq \cdots \leq \k_{n-1}(x)$ be the ordered principal 
curvatures of $\bo$, i.e., the eigenvalues of the second fundamental form w.r.t.\ the inner normal.
Then $\bo$ is said to be {\bf strictly $k$-convex} if $\k_k(x) >0$ for all $x\in \bo$.
For the following we will need  the boundary $\bo$ to be both strictly $k$-convex and strictly $(n-k+1)$-convex
by the duality in the arguments. Of course this condition for the lesser  principal curvature
is enough.

\Theorem{\DD.2. Part C} {\sl
Suppose that $\bo$ is strictly $\min\{k, n-k+1\}$-convex.  Let $\psi\in C(\ob)$ be an inhomogeneous
term with values in $I_k$, i.e.,
$$
\psi(\ob)\ \ss\ I_k \ =\ \left( (n-2k){\pi\over 2}, (n-2(k-1)){\pi\over 2}      \right).
$$
Then existence and uniqueness hold for the Dirichlet problem for any continuous boundary 
values $\vf\in C(\bo)$.
}

 This is the best result so far for the inhomogeneous
Dirichlet problem in the continuous setting.

It is an interesting question whether Theorem \DD.2 could be extended to
 general functions $\psi:\ob\to (-n{\pi\over 2}, n{\pi\over 2})$.

We would like to point out that solutions to (\II.5) have a very nice geometric 
interpretation which goes back to [CG].  
If $u$ is a smooth function on a domain  $\O\ss\rn$, then the graph of $Du$ 
on $\O\times \rn$ is a Lagrangian submanifold  $L$ of $\rn\times \rn$.  This gives us 
a {\bf phase function} $\theta: L \to \bbr/2\pi \bbz$ for the tangent planes of $L$
by setting ${\rm Re}\{dz_1\wedge\cdots\wedge dz_n\} = e^{ i \theta}$.
 Furthermore, as pointed out  in [CG], this phase function satisfies the equation
 $$
 \nabla \theta = - JH  \qquad{\rm on}\ \ L
 \eqno{(\II.6)}
$$
where $H$ is the mean curvature vector of $L$.  
This proves that a Lagrangian submanifold is minimal if and only if it has constant phase
(and is therefore special Lagrangian).
The equation (\II.6) was left as an exercise in [CG].
However, since this paper is devoted to the SL potential equation, we have inserted a 
proof  in Appendix A.

Now note that if $u$ is a smooth solution to (\II.5), then the phase function $\theta$ is just
the pull-back of $\psi$ to $L={\rm graph}(Du)$. In particular, that pull-back satisfies (\II.6).

The result in  [CPW], discussed above,    led us to show that the special Lagrangian potential operator satisfies
 the condition, in our paper  [IDP],  of being ``tamable''.  
 In Theorem \DD.1 we prove that $\tan\left\{ {1\over n} f(A)  \right\}$ is tame on $\bbf_\theta$
 = the inverse image under $f$ of any subinterval $[(n-2){\pi\over 2} + \d, n{\pi\over 2})$ of the top
 phase interval $I_1$. Here $\theta \equiv (n-2){\pi\over 2} +\d$ and $\d>0$.
  As a result we get a different proof
 of the result of [DDT] mentioned above and stated in Theorem \DD.2.
 
 With this same phase constraint, Ryosuke Takahashi [T, Thm.\  1.1], also applying
 the tangent function to $f(A)$, proved that the composition $\tan f(A)$ is concave on $\bbf_\theta$
 if $\theta > (n-2){\pi\over 2}$.  
 This important result allows Takahashi to apply Evans-Krylov with many  consequences.

 The SL subequation $\bbf_\theta \equiv \{f(A) \equiv \sum_i \arctan \l_i (A)\geq \theta\}$, 
 which is defined in terms of the eigenvalues of $A$,
is {\bf universal} in the sense that it defines many associated subequations as we now explain. 

\Def{\II.2} Let $\ggg:\Sym(\bbr^n)\to \bbr$ 
be a homogeneous polynomial of degree $m$ which is G\aa rding hyperbolic with respect to the
identity $I$. This means  that $\ggg(tI + A)$ is a   polynomial having all real roots  for every  $A\in \Symn$. 
The negatives of the roots are called the {\bf G\aa rding eigenvalues},  which we list in order
$$
\l^\ggg_1(A) \leq \cdots \leq \l^\ggg_m(A),
\eqno{(\II.7)}
$$
so  $\ggg(tI + A) = \ggg(I) \prod_{k=1}^m (t+\l^\ggg_k(A))$.  The open set $\G$
 where $\l^\ggg_1(A) >0$ (i.e., where all the G\aa rding eigenvalues are $>0$), is a convex cone containing
 $I$,  and is called the {\bf G\aa rding cone}.  It has the property that for each $k$
 $$ 
 \l^\ggg_k(A) \ <  \  \l^\ggg_k(A+B) \qquad \text{for all } \ B\in \G.
 \eqno{(\II.8)}
$$ ([HP, Thm. 2.1] or [HP2]).    In particular,  the {\bf G\aa rding branches}
 $$
 \bL^\ggg_k \ =\ \{ A : \l_k^\ggg(A) \geq 0\}
\eqno{(\II.9)}
$$  
satisfy $ \bL^\ggg_k  + \overline{\G} \ss  \bL^\ggg_k$ for all $k$.
If, in addition, $\cp \equiv \{A\geq 0\}\ss \overline{\G}$ ,we call $\ggg$ a {\bf G\aa rding-Dirichlet} polynomial.
In this case each of the G\aa rding branches is a subequation (see the beginning of the next section).

\Def{\II.3}    Given a G\aa rding-Dirichlet polynomial $\ggg$ we define the $\ggg$-{\bf special Lagrangian potential operator} (or $\ggg$-{\bf SL operator} for short)
$$
f^\ggg(A) \  \equiv \ \sum_{k=1}^m \arctan(\l^\ggg_k(A))
\eqno{(\II.10)}
$$
for $A \in \Symn$ with values in $(-m\pitwo, m\pitwo)$.  We define the subequation
$$
\bbf_\theta^\ggg \ \equiv\ \{A\in\Sym(\bbr^n) : f^\ggg(A) \geq \theta\}
 \quad \text{for each  \ } \theta \in  \left (-m\pitwo, m\pitwo \right ),
\eqno{(\II.11)}
$$
and the associated {\bf $\ggg$-special  Lagrangian potential equation}
$$
\partial \bbf_\theta^\ggg \ \equiv\ \{A\in\Sym(\bbr^n) : f^\ggg(A) = \theta\}.
\eqno{(\II.12)}
$$
That $\bbf_\theta^\ggg$ is a subequation follows from (\II.8) and the assumption $\cp\ss\G$.

The special phases and the phase intervals are defined exactly as above with $n$ replaced by
$m$.   The existence and uniqueness go through for these equations exactly as in the basic case.
(See Section \AA.)  Our boundary calculation also holds.

\noindent
{\bf GENERALIZED THEOREM \BB.1.}
 {\sl
Let $\ggg$ be a G\aa rding-Dirichlet polynomial as above.  Then the 
asymptotic subequation $\abbf_\theta^\ggg$ of\ \  $\bbf_\theta^\ggg$, for $\theta \in (-m{\pi  \over 2},  m{\pi  \over 2})$,
    is given as follows.

(1) \ \ If $\theta\in I_k$ ($k=1, ... ,m$), then
$$
\abbf_\theta^\ggg\ =\ \bL_k^\ggg.
$$

(2) \ \ If $\theta_k$ ($k=1, ... ,m-1$) is a special value, then}
$$
\abbf_{\theta_k}^\ggg\ =\ \bL_k^{\s_{m-1}^\ggg}.
$$

This result along with the analogue of  Theorem \BB.6 (Proposition \EE.3), will be proved in Section \EE.
However, the proofs are parallel to the basic case, and the reader should first study those arguments.

\noindent
{\bf Examples.}
On $\bbc^n$ we can take $\ggg$ to be the determinant of the complex symmetric part $\half(A-JAJ)$
of $A$, and we get the complex version of the special Lagrangian potential equation, 
which is related to mirror symmetry, as noted below. 
There is a similar quaternionic   case  with $\ggg$ equal to ${1\over 4} (A -IAI-JAJ-KAK)$. 
 In fact this process yields an infinite
array of equations, and is discussed in more detail in Section \EE.

 Another interesting  incarnation  of the SL potential equation is the following.  Let $u: \O \to \bbr$ 
 be a smooth function on a domain $\O\ss\rn$, and let $\Gamma \ss\O\times \bbr$ be its graph. 
 For each $x\in \O$, let $\kappa_1(x), ... , \kappa_n(x)$ be the eigenvalues of the second fundamental
 form, i.e., the {\sl principle curvatures},  of $\Gamma$ above $x$.  Then we can consider the subequation  
 ${\mathcal F}_{\theta, \rho}$  given  by  $\sum_{k=1}^n \arctan ( \rho \kappa_k) \geq \theta$
 for $\theta \in (-n {\pi\over 2}, n {\pi\over 2})$ and $\rho>0$.  This equation was first studied by Graham Smith  in [GS1].
 He considered the case where $\theta \in  [(n-1){\pi\over 2},  n{\pi\over 2})$ where  all  the  principal curvatures must be positive. The special case where  $\theta = (n-1){\pi\over 2}$ has a very nice geometric interpretation. 
 Actually, Smith considered this for hypersurfaces in a general riemannian manifold.  The paper [GS1] principally concerned rigidity and pre-compactness theorems.  However,  in [GS2] he considered the Dirichlet problem, for the same equations, on domains in Hadamard manifolds.  He went on to establish ``parametric'' versions of these theorems for  general convex curvature functions (see, for example, [GS3], [CS], [GS4]).

 Concerning the standard Dirichlet problem for this equation, even though it is constant coefficient
 (but not pure second order), comparison is an open question (although
a weak form of comparison does hold).  
Nevertheless, existence  was established in [DDR], for all $\theta$, 
on any domain in $\rn$  whose boundary satisfies the appropriate  ${\mathcal F}_\theta$  convexity.
Furthermore,  a version of our Theorem 3.1, explicitly computing this convexity condition, holds for this equation.
This is discussed in the last subsection of Section 6.

Using the results in [DDR] these results carry over to riemannian manifolds with
a topological $G$-structure.  This is discussed in Section \FF.

\noindent
{\bf Some Further Historical Notes.}  There is now a vast literature on the SL potential equation. 
The papers touch on many different topics, and we thought it might be interesting to
mention some examples.  However, this is  certainly not a complete history.

In Lemma 2.1  of  [Y2],  Yu Yuan proves that 
$$
\partial \bbf_\theta \ \text{ is convex (for some orientation) }
\qquad\iff\qquad 
|\theta|\geq (n-2){\pi\over 2}.
$$
For this reason $( (n-2){\pi\over 2}, n{\pi\over 2})$ is called the {\bf critical phase interval}.

 Dake Wang and Yu Yuan [WY1] showed that, for $m=2,3,...$,  there exist (viscosity) solutions $u_m$ to the SL Potential Equation  which are $C^{1, {1\over 2m-1}}$  but not $C^{1,\d}$ for $\d>{1\over 2m-1}$.  In fact these solutions are   analytic outside the origin. Graphing their gradients   gives special Lagrangian submanifolds
 with a isolated singularity at the origin.
 
 Recently,  J. Chen, R. Shankar and Y. Yuan proved that convex viscosity solutions to the 
 special Lagrangian potential equation are real
 analytic [CSY].

In another version of the boundary value problem,  Simon Brendle and Micah Warren [BW] proved
that if $\O_1, \O_2 \ss\ss\rn$ have smooth, strictly convex boundaries (second fundamental forms $>0$),
then there is a diffeomorphism $\Phi: \overline{\O_1}\to \overline{ \O_2}$ whose graph is special
Lagrangian (for some phase $\theta$).

There is a large literature concerning the SL potential equation and mean curvature flow.
A good survey is given by Wang [W].

People have worked on showing existence of Special Lagrangian submanifolds by 
minimizing volume among just among Lagrangians. (Recall from [CG] that Lagrangian
{\bf and} minimal implies special Lagrangian; see (\II.6) and Appendix A.)  This was started
by the work of Schoen and Wolfson [SW1,2].  It turns out to be quite subtle.  Wolfson [Wo] found an
example of a minimizer among Lagrangians   that was not minimal.  This lead to looking at mean curvature flow.
Here singularities   occur quite often in finite time (cf. A.  Neves [N1,3]).  See the surveys  [N2] and [W].

\def\pit{{\pi\over 2}}

There is a Bernstein-type theorem proved for $n=2$ by Lei Fu [F] and for general $n$ by Y. Yuan [Y1]
and also by Jost and Xin [JX].
It says that if $u:\rn \to \bbr$ is a global solution of the SL potential equation with phase
$|\theta| \in ((n-2)\pit, n\pit)$ (the critical phase interval), then $u$ is quadratic.

A degenerate form of the SL potential equation governs geodesics in the space
of positive graph Lagrangians in $\bbc^n$. The program for studying this space was
initiated by Jake Solomon, and he and Yanir Rubinstein [RS]  were able to solve this geodesic equation
with a continuous solution in the sense of [DD].
Matt Dellatorre expanded these results to manifolds with curvature $\leq 0$ [De].  It is very interesting 
that this program  is much like the program which lead to the solution of the Donaldson-Tian-Yau
conjecture [CDS]. Recently Collins and Yau [CY] have studied geodesics on an infinite dimensional
space which is mirror to Solomon's.  This is aimed at understanding the deformed Hermitian-Yang-Mills
equation which is ``mirror'' to the SL potential equation.

There is a fundamental article by N.  Hitchin [H] on the moduli space of special Lagrangians,
 and many important articles by  D. Joyce.  Some more recent articles are  [Jo1],  [Jo2] and [JLS].
Joyce's earlier work is surveyed in [Jo3] (together with a large overview of the field). 
 There is also a 
moment map point of view on this equation given by Donaldson [D].

The  SL potential  equation plays a big role in mirror symmetry.  
This began with the paper of A. Strominger, S.-T. Yau and E. Zaslow [SYZ]
which gave a very geometric picture of how mirror manifolds are connected.
Special Lagrangians and the SL potential equation are critical in this tableau. See,  for example, 
the very good articles [LYZ], [JY], [J],  [CJY], [CXY], [CSh].  The reader should consult these
sources, but a small insight comes from the following.  Let $(X,\o)$ be an $n$-dimensional K\"ahler manfold
and $a\in H^{1,1}(X, \bbr)$ a fixed (1,1)-homology class.  One is interested in finding
an element $\a\in a$ such that 
$$
{\rm Im}  \left ( e^{-i\theta} (\o+i\a)^n   \right) \ =\ 0.
$$
  The angle $\theta$ is determined topologically  by
  $$
\theta \ =\ {\rm arg} \left\{  \int_X (\o+i\a)^n  \right\}.
$$
This gives rise to a hermitian Yang-Mills equation 
$$
\Theta_\o(\a) = \sum_k \arctan(\l_k) \ \equiv \  \theta \ (  {\rm mod}\ 2\pi)
$$
where the $\l_k$'s are eigenvalues of an endomorphism $K:T^{1,0} X \to T^{1,0} X$
given by contracting by $\a$ and the dual of $\o$.  Of course the elements in $a$ all differ
from a given one $\a_0$ by $d d^c u$ for a function $u$ on $X$.

We point out that the work of R. Takahashi, referred to above,  has consequences for
some of these mirror symmetry results (see [T] for details).

Very recently Gao Chen confirmed a mirror version of the Thomas-Yau  conjecture (see [GC1-2]).

We want to thank the referee for his/her useful remarks, and particularly for noticing an oversight
in Proposition \BB.5.

%\vskip .3in
\vfill\eject

%%%%%%%%%%%%%%%%%%%%%%%%%%%%%%%%%%%%%%%%%%%%%%%%
%%%%%%%%%%%%%%%%%%%%%%%%%%%%%%%%%%%%%%%%%%%%%%%%
%%%%%%%%%%%%%%%%%%%%%%%%%%%%%%%%%%%%%%%%%%%%%%%%
%%%%%%%%%%%%%%%%%%%%%%%%%%%%%%%%%%%%%%%%%%%%%%%%
%%%%%%%%%%%%%%%%%%%%%%%%%%%%%%%%%%%%%%%%%%%%%%%%
%%%%%%%%%%%%%%%%%%%%%%%%%%%%%%%%%%%%%%%%%%%%%%%%
%%%%%%%%%%%%%%%%%%%%%%%%%%%%%%%%%%%%%%%%%%%%%%%%

\def\H{\bo}
\def\L{ \bL}

\noindent {\headfont \AB.  Geometric Conditions for Strict $\abbf_\theta$-Boundary Convexity}

The best boundary condition for existence in the $\bbf_\theta$-Dirichlet problem is the strict
 $\abbf_{|\theta|}$-convexity of $\bo$ at each point. As explained above, one needs in general 
 to verify convexity for both the subequation and its dual. However, the dual of $\bbf_\theta$ is $\bbf_{-\theta}$,
 and $\bbf_\theta \ss \bbf_{-\theta}$ if $\theta \geq 0$, so that one only need verify boundary convexity 
 for the smaller $\bbf_\theta$.
 
  The condition of strict $\abbf_{\theta}$-convexity can be computed entirely 
 in terms of the geometry of $\bo$.  Although we only need $\theta\geq 0$,
 we give the result for all $\theta$.

\Theorem{\AB.1}  {\sl 
 Fix $\theta \in (-n {\pi\over 2}, n{\pi\over 2})$.  Let $x\in \bo$ and let 
  $$
 \k_1\leq \cdots\leq \k_{n-1}
 $$
  be the principle curvatures of the second fundamental
 form  at $x$ w.r.t.\ the interior unit  normal $n$.  Then $\bo$ is strictly
$\abbf_\theta$-convex at $x$ if and only if:

(1) \ \ for $\theta \in I_k$ (the $k^{\rm th}$ phase interval) where $k< n-1$, one has
$$
\k_k\ >\ 0,
$$

(2) \ \  for $\theta = \theta_k = (n-2k){\pi\over 2}$, a special phase with  $1\leq k\leq n-2$, 
either

\hskip .3in (a) \ \ $\k_{k+1} > 0$ and $\k_k\geq 0$, or

\hskip .3in (b) \ \ $\k_{k+1} > 0$, $\k_k < 0$, and $\s_{n-1}(\k) \s_{n-2}(\k) <0$

\noindent
where $\s_\ell(\k)$ denotes the $\ell^{\rm th}$ elementary symmetric function of the $\k_i$'s.

(3)\ \  for $\theta = \theta_{n-1}$ or for $\theta \in I_n$, there is no condition.  (Every boundary is strictly 
 $\abbf_\theta$-convex.)}

Thus for the smallest phases $-n{\pi\over 2} < \theta \leq (-n+2){\pi\over 2}$, where $\oa\bbf_\theta$ is large,
there is no condition.  However, as $\theta$ increases,  the subequation $\oa\bbf_\theta$ decreases, culminating in the most 
stringent case  $\oa\bbf_\theta = \cp$ when $\theta \in I_1$.

\noindent
{\bf Proof.}
Consider the diagonal matrix
 $$
A\ =\ 
 \left(
 \begin{matrix}
\  \k_1   &  \    & \ & \           \\
\   & \k_{2} & \  & \       \\
\     & \ & \k_{3} &    \    \\
\ &  \ &  \ &  \cdots &      \\
\ &  \ &  \ &  \ & \k_{n-1}     \\
\ &\ &\ &  \ & \   & t      
 \end{matrix}
 \right).
 $$
 written with respect to an orthonormal basis of principal directions in $T_x(\bo)$ and $n$.
The condition is that $A$ should lie in $\Int \abbf_\theta$ for all $t>>0$ (cf. Cor. 5.4 in [DD]).
We use Theorem \BB.6.  If $\theta\in I_k$ for $k\leq n-1$, we have $\k_k>0$. 
This proves (1).   If $\theta \in I_n$, then we only
need $t>0$ which is always true.
 If $\theta = \theta_k$, then either (a) or (b)
must hold, for $k<n-1$.  When $k=n-1$, the $n^{\rm th}$ eigenvalue is $t$, which is $>>0$,
and the  $(n-1)^{\rm st}$ eigenvalue $\k_{n-1}$ is either $\geq 0$ or $<0$.  In this second case
all $\k_i$'s are $<0$ and $\s_{n-1}(A) \s_n(A)<0$ for all $t$ large (see (\AB.1) below.)  \qed

\noindent
{\bf Definition \AB.2.}
If the $k^{\rm th}$ ordered principal curvature at $x\in\bo$  satisfies
$
\k_k  > 0
$
we say  that $\bo$ is strictly {\bf $k$-convex} at $x$. This means that $\bo$ has at least $n-k$ strictly positive
principle curvatures at $x$. Notice in particular that 
strictly 1-convex means the hypersurface $\bo$ is strictly convex near $x$ in the usual sense.

\noindent
{\bf Examples of Strictly $\ell$-Convex Boundaries.}  Let $M\ss \rn$ be a smooth compact submanifold of 
dimension $\ell-1$  (codimension  $k+1$ with $k=n-\ell$).  For $\e>0$
let 
$$
S_\e(M) \equiv\  \{x : {\rm dist}(x,M)=\e\}.
$$
For $\e$ sufficiently small, $S_\e(M)$ is regular and diffeomorphic to the normal sphere bundle of $M$.
It bounds the domain $\O_\e(M) =  \{x : {\rm dist}(x,M)\leq \e\}$.  
For $\e>0$ even smaller $S_\e(M) = \partial \O_\e(M)$ has the property that 
at each point the second fundamental form has at least $k$ strictly positive eigenvalues
(coming from the normal spheres).
That is, the hypersurface $S_\e(M)$ is strictly $(n-k)=  \ell$ convex.
Thus there are many many strictly $\L_\ell$-convex boundaries in $\rn$, and in fact in any $n$-manifold.
Here is a picture of a strictly 2-convex surface in $\bbr^3$.

     \centerline{   %\hskip 1.3in
           \includegraphics[width=.3\textwidth, angle=0,origin=c]{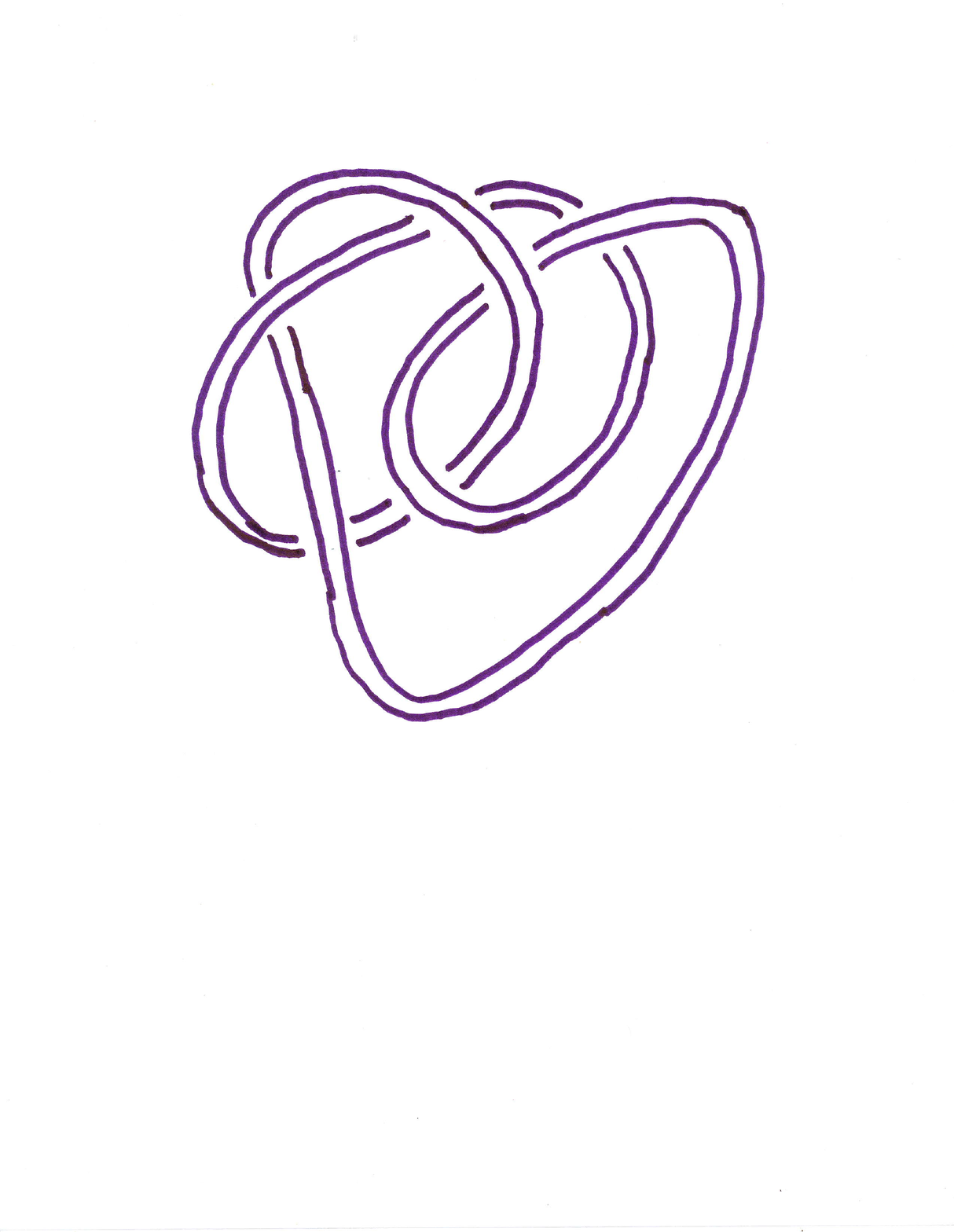}
          }
           \noindent  

\def\SFF{{\rm II}}
 
We now consider the case of a special phase $\theta = \theta_k$, $1\leq k\leq n-1$.
At a point  where (2)(a) holds, we have that $\bo$ is $k$-convex and also $\k_{k+1}>0$.

If  (2)(b) holds, i.e., if   $\k_{k+1} > 0$, $\k_k < 0$, and $\s_n(\k) \s_{n-1}(\k) <0$, we have
$$
{\s_{n-1}(A ) \over \s_n(A )} \ =\ {1\over t} + {1\over \k_{n-1}}+ {1\over \k_{n-2}} + \cdots + {1\over \k_{1}}.
\eqno{(\AB.1)}
$$
Now this is $<0$ for all large $t$ iff 
$$
{\s_{n-2}(\SFF) \over \s_{n-1}(\SFF)} \ =\   {1\over \k_{n-1}}+ {1\over \k_{n-2}} + \cdots + {1\over \k_{1}} 
\ < \ 0.
\eqno{(\AB.2)}
$$
 where
 $$
\SFF\ =\ 
 \left(
 \begin{matrix}
\  \k_{1}   &  \    & \ & \           \\
\   & \k_{2} & \  & \       \\
\     & \ & \k_{3} &    \    \\
\ &  \ &  \ &  \cdots & \   \\
\ &  \ &  \ &  \  & \k_{n-1}    
 \end{matrix}
 \right).
 $$
is the second fundamental form of $\H$.
\bigskip

\centerline
{\bf Examples for the  Special Phases.}

\noindent
{\bf Example \AB.3. ($k=1,  n=3$)}.  When is an oriented surface $\Sigma\ss\bbr^3$
strictly $\L_1^{\s_{2}}$-convex?  This happens exactly when either:

(1) \ \   The second fundamental form of $\Sigma$ is $\geq 0$, and so the Gauss curvature 
$K=\k_1\k_2\geq0$,  and the mean curvature $H=\k_1+\k_2>0$.   

(2)  Otherwise we have $K<0$ and $H > 0$.

\noindent
This last fact follows from (\AB.2) which says
$$
{\s_{1}(\SFF) \over \s_{2}(\SFF)} \ =\ {H\over K} \ <\ 0.
$$

\noindent
{\bf Example \AB.4.} ($k=1$ but general $n$)   The hypersurface $\H$ is 
strictly $\L_1^{\s_{n-1}}$-convex when either:

(1) \ \   The second fundamental form of $\H$ is $\geq 0$, and so the Gauss-Kronecker  curvature 
$K_{GK}\equdef \det(\SFF) = \s_{n-1}(\SFF)\geq0$,  and $\k_2>0$.

(2) \ \  Otherwise we have $K_{GK}<0$ and   $\s_{n-2}(\SFF) > 0$.

\noindent
{\bf Example \AB.5. ($k=2, n=4$)}.   When is a hypersurface $\H\ss\bbr^4$ strictly $\L_2^{\s_{3}}$-convex?
The boundary $\H$ has principle curvatures $\k_1\leq\k_2\leq\k_3$.  Then $\k_3>0$ 
and either $\k_2\geq0$ or $\k_2<0$.
 In this second case
 $$
 {\s_2(\SFF)\over \s_3(\SFF)} \ =\ {1\over \k_1}+{1\over \k_2}+{1\over \k_3} \ < \ 0.
 $$
Here the Gauss-Kronecker  curvature $\s_3(\SFF) >0$,  and $\s_2(\SFF) <0$.

As mentioned above, for the $\bbf_\theta$ Dirichlet Problem we need the boundary to be
strictly $\abbf_{|\theta|}$-convex at each point.  Here there is a difference between $n$ even and
$n$ odd.

\noindent
{\bf Suppose $n=2m$.}   Then we have
$$
I_k \ =\ ((m-k)\pi, (m-k+1)\pi)
$$
and the non-negative special values are $(m-k)\pi$ for $1\leq k\leq m$.  Note in particular that
$0$ is a special value.

\noindent
{\bf Suppose $n=2m+1$.}   Then we have
$$
I_k \ =\ ((m-k +\half)\pi, (m-k+\smfrac 3 2)\pi)
$$
and the non-negative special values are $(m-k)\pi$ for $1\leq k\leq m$.  

   \centerline{   %\hskip 1.3in
           \includegraphics[width=.9\textwidth, angle=270,origin=c]{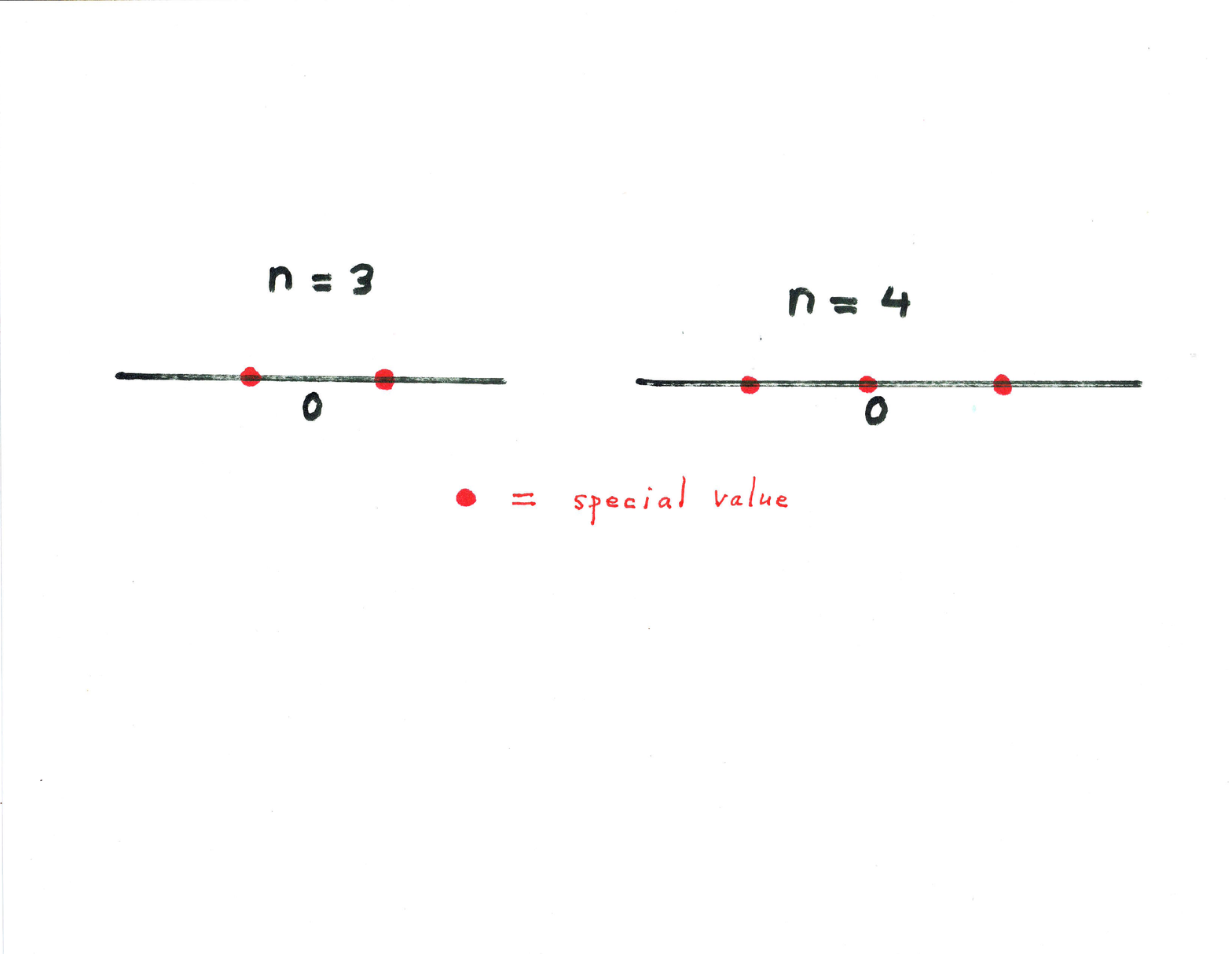}
          }
\noindent
{\bf Remark \AB.5. (Caffarelli - Nirenberg - Spruck).} In the ground-breaking paper [CNS]
the authors considered the case where 
$$
\theta \ \equiv \ 0 \, ({\rm mod} \ \pi)
$$
and solved the Dirichlet problem for $\theta$ in the uppermost interval.  For $n$ odd,
the boundary condition was strict convexity ($\k_1>0$).
However, for $n$ even, this $\theta$ is the largest special value.  From Theorem \BB.1 part (2)
the strict boundary convexity condition is either $\k_2>0$ and $\k_1\geq 0$ or $\k_2>0$, $\k_1 <0$ and
$\s_{n-1}\s_{n-2} <0$.  In the second case $\s_{n-1}<0$ and so $\s_{n-2}>0$. In the first case
we also have $\s_{n-2} >0$ because $\k_2>0$. Thus the strict convexity condition is exactly
$$
\s_{n-2}(\k_1, ... , \k_{n-1}) \ >\ 0.
$$
For a connected boundary $\bo$ of a compact domain, [CNS] had an elegant argument.
There is always a point $x\in\bo$ where the second fundamental form $\SFF$ is $>0$. Therefore 
$\SFF$ maps into the highest connected component of $\{\s_{n-2} \neq0\}$ (the one containing the identity).
Since $\s_{n-2}=0$ on the boundary of this region $b$ cannot touch that boundary by the strictness assumption.

\Remark {\AB.6}  Perhaps it is worth noting that the proof of Theorem \AB.1 can be interpreted as computing
the restriction of the subequation  $\abbf_{\theta}$ to the hyperplane $T_x\bo$.   
From this perspective, the condition in Theorem \AB.1 imposed on the second fundamental form
${\rm II}_{x, \bo}$ is that ${\rm II}_{x, \bo}$ must belong to the $\Int\left(\abbf_{\theta}\bigr|_{T_x\bo}   \right)$
of the restricted subequation $\abbf_{\theta}\bigr|_{T_x\bo} \ss\Sym(T_x \bo)$  on $T_x\bo$ (see [Rest]).
So, for example, in part (1) where  $\abbf_{\theta} = \L_k$,  we have computed that 
$\L_k \bigr|_{T_x\bo} = \L_k(T_k\bo)$,  the $k^{\rm th}$ branch of det = 0 on $T_x\bo$,
and the requirement is that ${\rm II}_{x, \bo} \in \Int \L_k(T_x\bo)$.

\vskip .3in
%\vfill\eject

%%%%%%%%%%%%%%%%%%%%%%%%%%%%%%%%%%%%%%%%%%%%%%%%
%%%%%%%%%%%%%%%%%%%%%%%%%%%%%%%%%%%%%%%%%%%%%%%%
%%%%%%%%%%%%%%%%%%%%%%%%%%%%%%%%%%%%%%%%%%%%%%%%
%%%%%%%%%%%%%%%%%%%%%%%%%%%%%%%%%%%%%%%%%%%%%%%%
%%%%%%%%%%%%%%%%%%%%%%%%%%%%%%%%%%%%%%%%%%%%%%%%
%%%%%%%%%%%%%%%%%%%%%%%%%%%%%%%%%%%%%%%%%%%%%%%%
%%%%%%%%%%%%%%%%%%%%%%%%%%%%%%%%%%%%%%%%%%%%%%%%

\noindent {\headfont \AA.   Preliminaries}

For the purposes of this paper a {\bf subequation} is a closed subset $\bbf \ss \Symn$
which is $\cp$-monotone, i.e., 
$$
\bbf + \cp \ \ss\ \bbf
$$
where $\cp \equiv \{A\geq 0\}$.    Associated to a subequation $\bbf$ there is a generalized potential theory,
which starts with the following concept. Given a domain $\O\ss\rn$ let $\USC(\O)$ be the set of functions
$u:\O\to [-\infty, \infty)$ which are upper semi-continuous. Then a function $u\in \USC(\O)$ is 
$\bbf$-{\bf subharmonic} if for each $x\in\O$ and each $C^2$-function $\phi\geq u$ in a neighborhood of
$x$ with $\phi(x)=u(x)$, we have $D_x^2 \phi\in\bbf$.

\vskip.3in

\centerline{\bf    The Special Lagrangian Potential Equation  }

The subequations of concern in this paper are defined as follows.

\Def{\AA.1}  The {\bf Special Lagrangian Potential Operator}, defined by
$$
f(A) \ \equiv\ \tr \left\{ \arctan A\right\}, \qquad {\rm for}\ \  A\in \Symn,
\eqno{(\AA.1a)}
$$
takes values in the interval $(-n \pitwo, n\pitwo)$.  Given a phase $\theta \in (-n \pitwo, n\pitwo)$, let 
$$
\Fth \ \equiv\ \{A\in \Symn : f(A)\geq \theta\}
\eqno{(\AA.1b)}
$$
denote the {\bf SL-potential subequation of phase $\theta$}.  The associated equation
$$
\partial \Fth \ \equiv\ \{A\in \Symn : f(A) = \theta \}
\eqno{(\AA.1c)}
$$
is called the  {\bf SL-potential equation of phase $\theta$}, and is usually written as
$$
\tr \{ \arctan \, D^2 u\}\ =\ \theta.
\eqno{(\AA.1d)}
$$

\vskip.3in

\centerline{\bf    The Asymptotic Interior }

There are two canonical subequations associated to each subequation $\bbf$:
the {\bf monotonicity subquation} ${\bf M}_{\bbf} \equiv \{A : \bbf+A\ss\bbf\}$, and the 
{\bf asymptotic subequation} $\abbf$, which determines the appropriate
 strict boundary convexity condition for existence in the Dirichlet problem.
 One of the main points of this paper is to compute $\aFth$ for the SL subequation $\Fth$.

 There are two approaches to defining and computing the asymptotic subequation $\abbf$
 for a  general subequation $F\ss \Symn$. Both approaches are enhancements of ideas   from [CNS].
 For the first approach we refer the reader to [DD, Sec. 5] where a rather complete description 
 of  various  cones which lead to $\abbf$ is given, along with examples showing the necessity 
 of the indirect approach using these auxiliary cones. (In [DD] $\abbf$ is called the {\sl asymptotic ray set} of $F$.)

Actually it is better to define and compute the   interior of $\abbf$
since it is really this which determines the strict boundary convexity.

This second approach was given in [DDR].  It has two advantages. First it is direct and does not require 
the auxiliary cones.  Secondly, it defines the interior of $\abbf$, which we provisionally denote by $\AI$.
Then $\abbf$ is defined to be its closure.  This is the approach taken here.

\Def{\AA.2}  Given a subequation $\bbf$ the {\bf asymptotic interior of} $\bbf$ is defined to be
$$
\AI \ \equiv\ \{A\in\Symn : \text{ $\exists\, \e>0$ and $t_0$ s.t.\ $t(A-\e I)\in \bbf \ \forall \, t\geq t_0\}$}.
%\eqno{(\AA.2a)}
$$
$
\text{The closure $\abbf$ of $\AI$ is called the {\bf asymptotic subequation} for $\bbf$.}
%\eqno{(\AA.2b)}
$

\Prop{\AA.3} {\sl
Let  $\bbf$ be a subquation which is $\neq \emptyset$ and not equal to $\Symn$.
Then $\AI$ is an open cone with vertex at the origin, and it is $\cp$-monotone.
Its closure is a subequation $\abbf$ also $\neq \emptyset$ or $\Symn$.
Finally, the interior of $\abbf$ is equal to the provisional  set $\AI$ (i.e., $\Int \abbf = \AI$).
}

\pf
It is easy to see from its definition that $\AI$ is a cone with vertex at the origin and that it is $\cp$-monotone.
Let
$$
{\cal N}_\e (A) \ \equiv\ A-\e I + \Int \cp
\eqno{(\AA.2)}
$$
denote the ``$\e$-neighborhood of $A$''.  Then the definition of $\AI$ can be reformulated as
$$
\AI \ =\ \{A :   \exists\  \e>0\ {\rm and}\  t_0 \ {\rm s.t.}\ t{\cal N}_\e (A) \ss\bbf \ \forall\, t\geq t_0\}.
\eqno{(\AA.3)}
$$

Now if $A\in\AI$, that is, $t {\cal N}_\e\ss F$ for all $t\geq t_0$, then 
$\forall\, B\in {\cal N}_{\e\over 2}(A)$, we have $t {\cal N}_{\e\over 2}(B)
\ss t{\cal N}_{\e}(A) \ss F$ for all $t\geq t_0$, which proves 
 that $\AI$ is open.  If $-I\in\AI$, then from (\AA.3)
 we see that  $\bbf = \Symn$ contrary to assumption.
 
 Since $\G \equiv \AI$ is an open $\cp$-monotone set which is non-empty and not
 equal to $\Symn$,  the remainder of the proposition is a consequence of the following
 more general fact.
 
 \Lemma{\AA.4}  {\sl 
 Suppose $\G$ is any open subset of $\Symn$ which is $\cp$-monotone, i.e., $\G+\cp\ss\G$.
 Then there exists a unique subequation $F\ss \Symn$ with $\G=\Int F$, or otherwise said, 
 $F=\overline \G$ is a subquation, and $\G$ is its interior.  Moreover, if $\G \neq \emptyset, \Symn$, then $F\neq \emptyset,  \Symn$.
 }
 
 \pf
 If $A\in F\equiv \overline\G$ and $P\geq 0$, then $A=\lim_{j\to \infty} A_j, \ A_j\in \G$ and hence $A+P = 
    \lim_{j\to \infty} (A_j +P)$ with $A_j+P \in \G$, that is, the closure $F$ of $\G$ is $\cp$-monotone, i.e., 
    $F$ is  a   subequation.
    
    Now $\G$ is an open set $\ss F$, hence $\G\ss \Int F$.  Finally, suppose $A\in \Int F$.  If $A\notin \G$,
    then $(A -  \Int\cp)\cap \G=\emptyset$.  Otherwise, $A-P\in\G$, $P > 0 \Rightarrow A\in \G$. 
    Since $A-\Int\cp$ and $\G$ are open disjoint sets, $A-\Int \cp$ and $F=\overline \G$ are disjoint.
    However, $A\in \Int F \Rightarrow A-\e I \in F$ for $\e>0$ small, which contradicts this intersections being empty.
    \qed

\vskip.3in

\centerline{\bf    Branches}

Let $\l_1(A)\leq \cdots \leq \l_n(A)$, or simply $\l_1\leq \cdots \leq \l_n$, denote the ordered eigenvalues
of $A\in\Symn$.  Recall that the $k^{\rm th}$ branch of the Monge-Amp\`ere equation $\det (A)=0$
is the subequation $\bL_k$ defined by
$$
\bL_k\ \equiv\ \{A : \l_k(A) \geq 0\}\qquad k=1, ... ,n.
\eqno{(\AA.4)}
$$

This construction extends to general equations $\ggg(A)=0$ where $\ggg(A)$ is a degree $m$ homogeneous 
polynomial on $\Symn$ which is G\aa rding hyperbolic in the direction of the identity $I$
(that is, all roots of $\ggg(tI+A)$ are real for $A\in \Symn$).  
The negatives of these roots are called the {\sl G\aa rding $\ggg$-eigenvalues}. 
%By renormalizing we can assume $p(I)=1$.
The ordered $\ggg$-eigenvalues  $\l_1^\ggg( A)\leq \cdots \leq\l_m^\ggg(A)$   determine the 
$k^{\rm th}$ branches $\bL_k^\ggg$ of the equation $\ggg(A)=0$ by 
$$
\bL_k^\ggg\ \equiv\ \{A : \l_k^\ggg(A) \geq 0\}\qquad k=1, ... ,m.
\eqno{(\AA.5)}
$$
The first (or principle) branch is the closure of the open convex cone $\G$  defined by 
$\l_k^\ggg(A) > 0$ for all $k$.  This open cone, $\G = \Int\bL_1^\ggg$, is called the 
{\bf G\aa rding cone} for $\ggg$, and $\overline{\G} = \bL_1^\ggg$ is the 
{\bf closed G\aa rding cone} for $\ggg$.

Two cases of this construction are pertinent to understanding  strict boundary convexity for 
the Special Langrangian potential equation.  The first is the branches
$$
\bL_1 \ \ss\ \bL_2 \ \ss\   \cdots \ss \ \bL_n   
\eqno{(\AA.6)}
$$
of the Monge-Amp\`ere equation $\ggg(A) \equiv \det(A) = \s_n(A)=0$.   
Here the $\ggg$-eigenvalues of $A$ %$\l^{\det}(A)$, 
are just the standard eigenvalues  $\l_1, ... , \l_n$ of $A$.
The second is the set of branches
$$
    \bL_1^{\s_{n-1}} \ \ss\    \bL_2^{\s_{n-1}}  \ \ss\   \cdots \ss \    \bL_{n-1}^{\s_{n-1}} 
\eqno{(\AA.7)}
$$
of the equation $\ggg(A) \equiv \s_{n-1}(A)=0$ (where $\s_k$ is the $k^{\rm th}$ elementary symmetric 
function of the eigenvalues of $A$).

The second case is complicated by the fact that the $\s_{n-1}$-eigenvalues, 
$\l_k^{\s_{n-1}}( A)$,
cannot be computed in 
terms of the standard $\s_n$-eigenvalues $\l_1,...,\l_n$ of $A$.
 Note  that for any hyperbolic polynomial $\ggg$, after renormalization, we have
$$
\ggg(tI+A) \ =\ t^m + \s_1^\ggg(A) t^{m-1} + \cdots +  \s_{m-1}^\ggg(A) t + \s_m^\ggg(A)
$$
which  defines the $\s_k^\ggg(A)$'s.
As a result  we have 
$$
\s_{m-1}^\ggg(A) \ =\ {1\over m} {d\over dt}\ggg(tI+A)\bigr|_{t=0},
$$
so the eigenvalues of  $\s_{m-1}^\ggg(A)$ are  the critical points of $\ggg(tI+A)$ which is  $\det(tI+A)$ in our case.

Although  the $\s_{n-1}$-eigenvalues, 
$\l_k^{\s_{n-1}}( A)$,
cannot be computed in 
terms of the standard eigenvalues of $A$, 
 the condition $\l_k^{\s_{n-1}}( A)  \geq 0$ can be described by inequalities involving these
eigenvalues.  This  key result, Proposition \BB.5, was stated in the introduction and 
 is proved near the end of  Section \BB.

\vskip.3in

\centerline{\bf   Some General Results on the Pure Second-Order Dirichlet Problem
}

Here we recall  what is known regarding the Dirichlet problem for an arbitrary
subequation $\bbf\ss \Symn$  (defined at the beginning of  this section).  
The results all apply to the equation $\bbf_\theta$.
However, in Section 4 we get finer results for $\bbf_\theta$ by using Theorem \BB.1 to implement
the general results here.

Suppose $\O$ is a bounded domain in $\rn$, and   $X\ss\rn$ is an arbitrary   open subset.

\Def{\AA.5} A function $h\in C(X)$ is {\bf  $\partial\bbf$-harmonic}
if $h$ is $\bbf$-subharmonic and $-h$ is $\wt\bbf$-subharmonic where
$\wt\bbf = \ \sim(-\Int \bbf)$ is the dual subequation. (Note that $\partial \bbf = \bbf \cap (-\wt\bbf)$.)

\Def{\AA.6}   We say that 
{\bf existence for the (DP) for $\partial \bbf$ holds on $\O$} if for all prescribed boundary functions
$\vf \in C(\bo)$ there exists $h\in C(\ob)$ satisfying

 (a) \ \ $h\bigr|_{\O}$ is $\partial\bbf$-harmonic, and

 (b) \ \ $h\bigr|_{\bo} \ =\ \vf$.

\Def{\AA.7}   We say that 
{\bf uniqueness for the (DP) for $\partial \bbf$ holds on $\O$} if for all prescribed boundary functions
$\vf \in C(\bo)$ there exists at most one  $h\in C(\ob)$ satisfying (a) and (b).

Consider also the following stronger versions of existence and uniqueness.

\Def{\AA.8}  If for all $\vf \in C(\ob)$ the Perron function %$h(x) \equiv \sup\{u(x) : u\in \USC(\ob),\
%u \ {\rm is } \ \bbf-{\rm subharmonic,\ and} \ u\bigr|_{\bo} \leq \vf\}$
$$
\text
{
 $h(x) \equiv \sup\{u(x) : u\in \USC(\ob),\
u$ is $\bbf$-{subharmonic, and}  $u\bigr|_{\bo} \leq \vf\}$
}
$$
is in $C(\ob)$ and satisfies (a) and (b), then we say that
{\bf Perron existence holds for $\partial \bbf$ on $\O$}.

\Def{\AA.9} If for all $u,v\in \USC(\ob)$ with  $u\bigr|_{\O}$ $\bbf$-subharmonic and
 $v\bigr|_{\O}$ $\wt\bbf$-subharmonic,
 $$
 u+v \ \leq\ 0 \quad {\rm on} \ \bo
 \qquad\Rightarrow\qquad
  u+v \ \leq\ 0 \quad {\rm on} \ \ob,
 $$
 then we say that {\bf comparison holds for $\bbf$ on $\O$.}
 
 Obviously Perron existence implies existence, and comparison implies uniqueness.
 
In [DD] comparison was shown to always hold.

\Theorem{\AA.10. (Comparison)}  {\sl
For all bounded domains $\O\ss\rn$ and $u,v\in \USC(\ob)$  with $u\bigr|_{\O}$ $\bbf$-subharmonic and
 $v\bigr|_{\O}$ $\wt\bbf$-subharmonic,
  $$
 u+v \ \leq\ 0 \quad {\rm on} \ \bo
 \qquad\Rightarrow\qquad
  u+v \ \leq\ 0 \quad {\rm on} \ \ob.
 $$
}

(See Remark 4.9 in [DD] for the proof that the definitions of $\bbf$ and $\wt\bbf$-subharmonicity
agree with the appropriate viscosity definitions.)

Theorem \AA.10 applies to the SL-potential equation $\bbf_\theta$ for arbitrary phase $\theta$
(see (\AA.1b)), extending the result of [CNS] for $(n-1){\pi\over 2} < \theta < n {\pi\over 2}$,
and settling the comparison/uniqueness question in the affirmative for all bounded domains.

A second proof of Theorem \AA.10 was given later in [DDR].

This leaves the existence question for $\bbf$, which is also covered by the following result of [DD].

\Theorem{\AA.11. (Perron Existence)}  {\sl
Suppose that $\O$ is a bounded domain with smooth boundary $\bo$.
If, at each point of $\bo$, the boundary is both
$$
\text
{
$\bbf$-strictly convex, and $\wt\bbf$-strictly convex,
}
\eqno{(\AA.8)}
$$
then Perron existence for $\partial \bbf$ holds on $\O$.
}

The definition of strict convexity in (\AA.8) has many equivalent formulations.
The key is the asymptotic interior $\Int \oa \bbf$
 of $\bbf$ (see Definition \AA.2).

Here are two of the equivalent definitions of strict $\bbf$-convexity for $\bo$.
Let ${\rm II}_{x, \bo} $ denote the second fundamental form of $\bo$ with respect
to the interior unit normal $n=n_x$, and let $P_n$ be orthogonal projection onto the line $\bbr\cdot n$.
Then for each $x\in \bo$, {\bf strict $\bbf$-convexity} is the requirement that 
$$
{\rm II}_{x, \bo} + tP_{n} \ \in \ \Int \oa \bbf
\ \ {\rm for} \ t \gg 0,
\eqno{(\AA.9a)}
$$
This definition is equivalent to the following.  
$$
\begin{aligned}
&\text
{
$\exists\, \rho\in C^\infty(\ob)$ with $\O=\{\rho<0\}$,  and  $\rho=0$ and $\nabla \rho\neq 0$ on $\bo$,
}
\\
&\qquad\qquad\qquad \text
{
such that \ \ \  $D^2_x \rho \in \Int \oa\bbf\ \forall x\in \ob$.
}
\end{aligned}
\eqno{(\AA.9b)}
$$
See [DD] Lemma 5.3 and Corollary 5.4 for the equivalence of (\AA.9a) and (\AA.9b).

Summarizing by combining Theorems \AA.10 and \AA.11, we have for a general 
subequation $\bbf\ss\Symn$ the following.

\Theorem{\AA.12. (Dirichlet Problem)}  {\sl
Suppose that $\O$ is a bounded domain with $C^\infty$ boundary $\bo$ with is
both $\bbf$ and $\wt\bbf$ strictly convex.  Then both Perron existence and comparison
hold for $\bbf$ on $\O$.  In particular, both existence and uniqueness hold for the $\bbf$
Dirichlet problem on $\O$
}

In particular all this applies to the SL potential equation of arbitrary phase.
Now we come to the main new result of this paper.

\vskip .3in
%\vfill\eject

%%%%%%%%%%%%%%%%%%%%%%%%%%%%%%%%%%%%%%%%%%%%%%%%
%%%%%%%%%%%%%%%%%%%%%%%%%%%%%%%%%%%%%%%%%%%%%%%%
%%%%%%%%%%%%%%%%%%%%%%%%%%%%%%%%%%%%%%%%%%%%%%%%
%%%%%%%%%%%%%%%%%%%%%%%%%%%%%%%%%%%%%%%%%%%%%%%%
%%%%%%%%%%%%%%%%%%%%%%%%%%%%%%%%%%%%%%%%%%%%%%%%
%%%%%%%%%%%%%%%%%%%%%%%%%%%%%%%%%%%%%%%%%%%%%%%%
%%%%%%%%%%%%%%%%%%%%%%%%%%%%%%%%%%%%%%%%%%%%%%%%

\noindent {\headfont \BB.  Computing the Asymptotic Interior of $\bbf_\theta$}     

\vskip.2in

\centerline{\bf Phases (or Values) of the Special Langrangian Potential Operator}

We now consider the operator $f(A) \equiv \tr \{ \arctan (A)\}$ defined on $\bbf \equiv \Symn$.  
This operator $f$ has values precisely in the open interval $(-n{\pi  \over 2},  n{\pi  \over 2})$.
There are $n-1$ 
$$
{\bf Special\ Phases (or Values):}   \hskip .3in  \theta_k \ =\ (n-2k) {\pi\over 2}, \quad k=1,..., n-1.  
\eqno{(\BB.1)}
$$
Removing these $n-1$ special values, the remaining set of values is 
the disjoint union of $n$ open
$$
{\bf  Phase \ Intervals:} \ \ 
I_k \ =\  \left ( (n-2k) {\pi\over 2}, \ (n-2(k-1)) {\pi\over 2} \right), \ k=1,...,n.
\eqno{(\BB.2)}
$$

\Theorem {\BB.1}  {\sl
The asymptotic subequation $\abbf_\theta$ of\ \  $\bbf_\theta$, for $\theta \in (-n{\pi  \over 2},  n{\pi  \over 2})$,
    is given as follows.

(1) \ \ If $\theta\in I_k$ ($k=1, ... ,n$), then
$$
\abbf_\theta\ =\ \bL_k.
$$

(2) \ \ If $\theta_k$ ($k=1, ... ,n-1$) is a special value, then}
$$
\abbf_{\theta_k}\ =\ \bL_k^{\s_{n-1}}.
$$

For the proof of both parts we will need  the asymptotic expansion for $f(tA)$ as $t\to \infty$.
For the proof of part (2) we will also need Proposition \BB.5 below describing $\bL_k^{\s_{n-1}}$.
Define
$$
q(A) \ =\ \text{the number of strictly negative eigenvalues of } \ A.
\eqno{(\BB.3)}
$$

\Lemma{\BB.2}  {\sl
Suppose $A$ is non-degenerate ($\s_n(A)\neq 0$).  Then}
$$
\begin{aligned}
f(tA) \ &=\ \tr\{\arctan(tA)\}\   \\
&=\  (n-2q(A) )   \pitwo \ -\ {1\over t} {\s_{n-1}(A) \over \s_n(A)} + O\left( {1\over t^3}\right) \quad {\rm as} \ \ t\to\infty.
\end{aligned}
\eqno{(\BB.4)}
$$

\Cor{\BB.3} {\sl
Suppose $\s_n(A) \neq 0$.  Then

(1)
$$
\lim_{t\to\infty} f(tA) \ =\ (n-2q(A))  \pitwo,
\eqno{(\BB.5)}
$$
and if also $\s_{n-1}(A) \neq 0$, then 

(2) \ \ $f(tA)$ either strictly decreases to $ (n-2q(A))  \pitwo$ or 
strictly increases to $ (n-2q(A))  \pitwo$  depending on whether
$\s_{n-1}(A)$ and $\s_n(A)$ have opposite signs or the same sign respectively.
}

\noindent
{\bf Proof of Lemma \BB.2.}   The asymptotic expansion for the $\arctan t$ as $t\to\infty$ is 
$$
\arctan(t) \ = {\pi\over 2} - {1\over t} + {1\over 3 t^3}  + \cdots  \quad {\rm as} \ \ t\to\infty.
\eqno{(\BB.6)}
$$
(Note that $\arctan(t)$ and $-\arctan(1/t)$ have the same derivative so that 
$\arctan(t) = \pi/2 -\arctan(1/t)$, and  for $s>0$ small, $\arctan(s) = s -{1\over 3} s^3 +\cdots$.)
The number of strictly positive eigenvalues of $A$ is $n-q(A)$ since $A$ has no zero eigenvalues.
Therefore,
$$
\begin{aligned}
f(tA) \ &=\  \tr\{\arctan(tA)\} \ =\ \sum_{k=1}^n \tr\{\arctan(t \l_k)\}  \\
&=\   (n-2q(A))  \pitwo \ -\ {1\over t} \left(    \sum_{k=1}^n {1\over \l_k}  \right) + O\left( {1\over t^3}\right).
\end{aligned}
$$
Since 
$$
\sum_{k=1}^n {1\over \l_k} \ =\ {\s_{n-1}(A) \over \s_n(A)},
$$
this completes the proof.\qed

Because of Definition \AA.2 and Proposition \AA.3, by taking closures/interiors we have the following equivalent version of Theorem \BB.1.

\noindent
{\bf   Theorem \BB.1$'$.}    

(1)$'$ \ \ If $\theta\in I_k$ ($k=1, ... ,n$), then
$$
\AIth\ =\  \Int \bL_k.
$$

(2)$'$ \ \ If $\theta_k$ ($k=1, ... ,n-1$) is a special value, then
$$
{\oa{\Int} \, \bbf}_{\theta_k}\ =\  \Int \bL_k^{\s_{n-1}}.
$$

We will make use of the following equivalent ways of defining the $k^{\rm th}$ branch of 
the equation $\det A=0$.

\Lemma {\BB.4} {\sl
One has that $A\in \bL_k \iff \l_k(A) \geq 0  \iff q(A) \leq k-1 \iff q(A) < k \iff
 (n-2q(A)) {\pi\over 2} \geq (n-2(k-1)) {\pi\over 2}  \iff 
 (n-2q(A)) {\pi\over 2} > (n-2k) {\pi\over 2}$.  One also has: 
 $A \in \bL_{k+1} \sim \bL_k \iff k\leq q(A) \leq k$, i.e., $q(A)=k$.
}

Next note that for each $A\in\Symn$, $\s_n(A-\e I)$ and $\s_{n-1}(A-\e I)$
are polynomials in $\e$, so that 
$$
\{ \e>0 : \s_n(A-\e I)=0\} \ \ {\rm and}\ \ \{ \e>0 : \s_{n-1}(A-\e I)=0\} \ \ {\rm are \  finite}.
\eqno{(\BB.7)}
$$

\noindent
{\bf Proof of Theorem \BB.1$'$ (1)$'$.}
Now assume $\theta \in I_k$, i.e., $(n-2k){\pi\over 2} < \theta < (n-2(k-1))\pitwo$.

If $A\in \Int \bL_k$, then $A_\e \equiv A-\e I \in \bL_k$ for $\e>0$ small.
By Lemma \BB.4, we have $A_\e \in \bL_k \iff (n-2(k-1)) \pitwo \leq (n-2q(A_\e))\pitwo$.  
Hence, $\theta < (n-2q(A_\e))\pitwo$.  We can assume that $A_\e$ is non-degenerate by (\BB.7),
so we can apply Corollary \BB.3 to obtain
$$
\theta \ <\ \lim_{t\to \infty} f(t A_\e).
$$
In particular, there exists $t_0$ such that $f(t(A-\e I))\geq \theta$ for all $t\geq t_0$ or equivalently
$t(A-\e I) \in \bbf_\theta$  for all $t\geq t_0$.
By Definition (\AA.2a) this proves $A\in \AIth$.

Conversely, suppose  $A\in \AIth$.  Then by  Definition (\AA.2a), there exists $\e>0$ and $t_0$ such
that $f(t(A-\e I))\geq \theta$ for all $t\geq t_0$, but $\theta \in I_k \Rightarrow (n-2k)\pitwo <\theta$.
By (\BB.7) we can assume $\s_n(A-\e I)\neq 0$, so Corollary \BB.3(1) applies to yield
$$
(n-2k)\pitwo \  < \ \theta \ \leq \ \lim_{t\to \infty} f(t A_\e) \ =\ (n-2q(A_\e))\pitwo.
$$
By Lemma \BB.4 this implies $A_\e = A-\e I \in\bL_k$.  Hence $A\in \Int \bL_k$.
This proves (1)$'$ and hence (1). \qed

To prove Part (2)$'$ we use Proposition \BB.5 below.

%%%%%%%%%%%%%%%%%%%%%%%%%%%%%%%%%%%%%%%%%%%%%%%%
%%%%%%%%%%%%%%%%%%%%%%%%%%%%%%%%%%%%%%%%%%%%%%%%
%%%%%%%%%%%%%%%%%%%%%%%%%%%%%%%%%%%%%%%%%%%%%%%%
%%%%%%%%%%%%%%%%%%%%%%%%%%%%%%%%%%%%%%%%%%%%%%%%
%%%%%%%%%%%%%%%%%%%%%%%%%%%%%%%%%%%%%%%%%%%%%%%%

\bigskip

\def\L{\bL}

\centerline{\bf Describing the Branches $\L_k$ and $\L_k^{\s_{n-1}}$  in Terms of }
\centerline{\bf Roots and Critical Points}

The $k^{\rm th}$ branch $\L_k$ of $\det(A)$ was defined in (\AA.4) as the set $\{\l_k(A)\geq 0\}$, $k=1, ..., n$.
The negatives of the eigenvalues  of $A$ will be referred to as the {\sl roots} of $A$.
More precisely, we adopt the following definitions.
$$
\begin{aligned}
&\text{
Let  \ \ \ $\ell\ \equiv\  n-k$ \qquad for  $k=1, ... , n$  \ \   and then   
}
\\
&\text{
let \ \ \ $r_{\ell+1} \equdef -\l_k,$ \qquad for $\ell+1=1,..., n$.
}
\end{aligned}
\eqno{(\BB.8)}
$$
Consequently, 
$$
r_1 \leq \cdots \leq r_{\ell+1} \leq \cdots \leq r_n,\qquad  \ell=0, ... , n-1.
$$
are the {\bf ordered roots}.  With this notation
$$
\L_k \ =\ \{A : \l_k(A) \geq 0\} \ =\ \{A : r_{\ell+1}(A) \leq 0\}\qquad\text{where $k+\ell=n$}
\eqno{(\BB.9a)}
$$
and
$$
\L_{k+1} \ =\ \{A : \l_{k+1}(A) \geq 0\} \ =\ \{A : r_{\ell}(A) \leq 0\}\qquad\text{where $k+\ell=n$}
\eqno{(\BB.9b)}
$$
In addition, their interiors are given by the corresponding strict inequalities, and their 
boundaries by equality.

The roots $r_j(A)$ are the roots of the monic polynomial $p_A(t)$ ( or simply $p(t)$) defined by
$$
\begin{aligned}
p(t) \ =\ \det(t I +A)  \ &=\ t^n +\s_1(A) t^{n-1} +\cdots +\s_n(A)  \\
 & =\ \prod_{k=1}^n(t+\l_k) \ =\ \prod_{\ell=1}^n(t-r_\ell)
 \end{aligned}
\eqno{(\BB.10)}
$$
As noted in Section \AA
$$
\s_{n-1}(tI+A) \ =\ {1\over n}{d\over dt} \det(tI+A) \ =\ {1\over n} p'(t).
$$
Hence the roots of $A$ for $\s_{n-1}$ are the {\sl critical points} of $p(t) = \det(tI+A)$,
which we denote by $c_1(A)\leq \cdots\leq c_{\ell}(A)\leq \cdots\leq c_{n-1}(A)$.
(It will not be necessary to label their negatives which are the $\s_{n-1}$-eigenvalues of $A$.)

The same reasoning used to prove (\BB.9) provides a description of the 
$k^{\rm th}$ branch of $\s_{n-1}=0$ (now the degree is $n-1$ instead of  $n$) in terms of the ordered critical points.
$$
\L_k^{\s_{n-1}} \ =\ \{A :  c_\ell(A) \leq 0\}  \qquad k+\ell = n\ \ {\rm and}\ \ k=1,..., n-1.
\eqno{(\BB.11)}
$$
Again,strict inequality defines the interiors and equality defines the boundaries.

With the translations above from eigenvalues to roots and critical points, 
we will make use of the following three elementary facts.

\noindent
{\bf Fact 1.}
$$
r_1  \leq c_1 \leq r_2 \ \leq \ \cdots \ \leq r_\ell \leq c_\ell \leq\ r_{\ell+1} \ \leq \cdots\ \leq \ c_{n-1}  \leq\ r_n.
$$
This fact implies the nesting
$$
\L_1 \ss \L_1^{\s_{n-1}} \ss \cdots \L_k \ss  \L_k^{\s_{n-1}} \ss \L_{k+1}
\ss \cdots \ss \L_n,\quad 1\leq k \leq  n-1.
\eqno{(\BB.12)}
$$

The second elementary fact eliminates all but two of the possibilities for  $r_\ell \leq c_\ell\leq r_{\ell+1}$.

\noindent
{\bf Fact 2.}  Either $r_\ell < c_\ell < r_{\ell+1}$  or   $r_\ell = c_\ell =  r_{\ell+1}$, \ \  for $\ell=1, ..., n-1.$

The third elementary fact is crucial.

\noindent
{\bf Fact 3.}  Assume $r_\ell < c_\ell < r_{\ell+1}$,  Then
$$
p(t)p'(t)  > 0 \ \ {\rm on} \ \ (r_\ell, c_\ell)
\and
p(t)p'(t)  <0 \ \ {\rm on} \ \ (c_\ell, r_{\ell+1}).
$$

\noindent{\bf Proof of Fact 3.}
Note that $p(t) \neq 0$ on $(r_\ell, r_{\ell+1})$ and $p'(t)\neq 0$ on both $(r_\ell, c_\ell)$ and $(c_\ell, r_{\ell+1})$.
In particular, $p(c_\ell)\neq 0$.  If $p(c_\ell)>0$, then $p(t)$ increases from 0 to $p(c_\ell)$ on $(r_\ell, c_\ell)$,
and hence $p(t)p'(t) >0$ on $(r_\ell, c_\ell)$.
 If $p(c_\ell)<0$, then $p(t)$ decreases from 0 to $p(c_\ell)$ on $(r_\ell, c_\ell)$,
and hence again $p(t)p'(t) >0$ on $(r_\ell, c_\ell)$.
The proof that $p(t)p'(t)<0$ on $(c_\ell, r_{\ell+1})$ is similar. \qed

     \centerline{   %\hskip 1.3in
           \includegraphics[width=.3\textwidth, angle=270,origin=c]{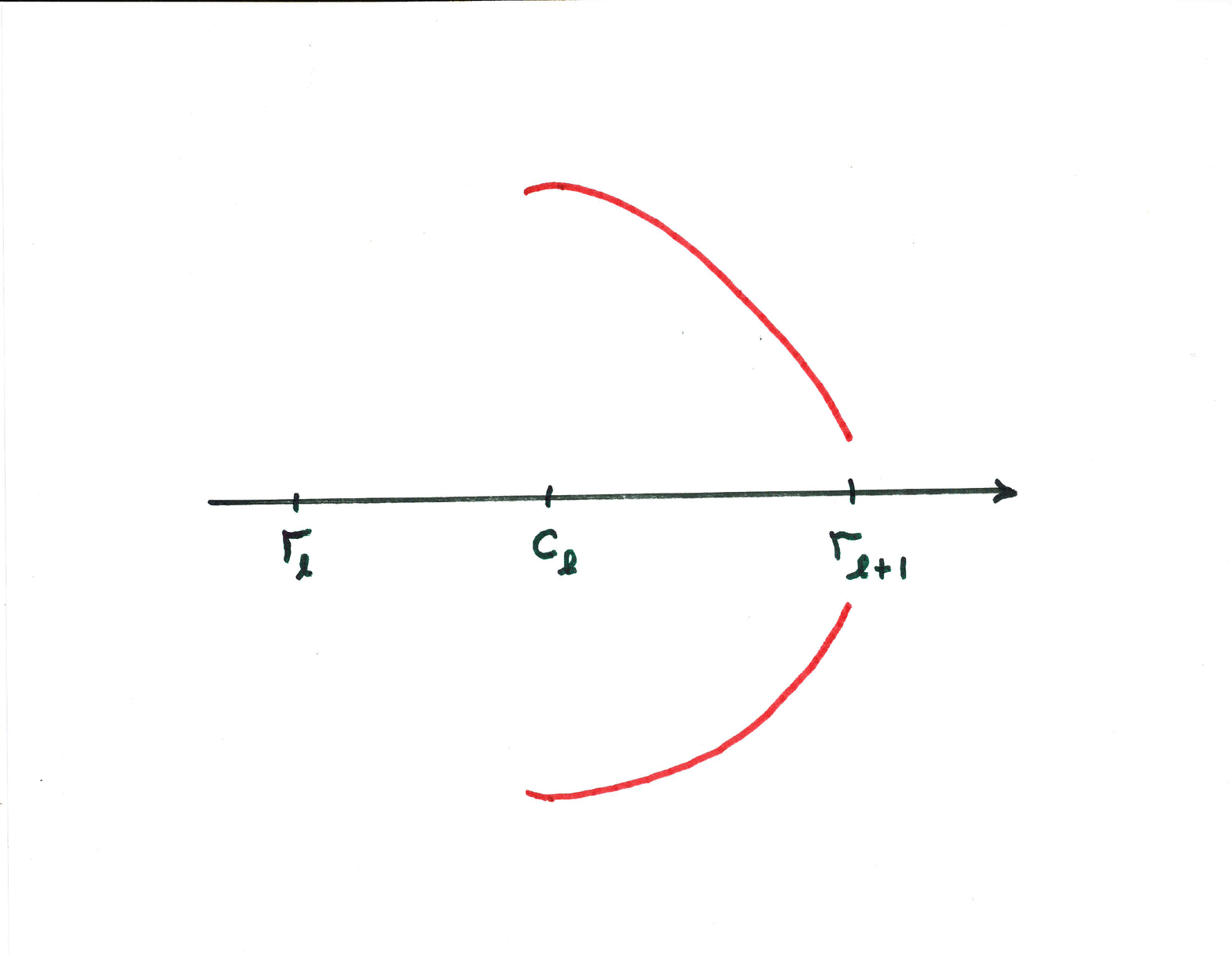}
          }

Each of the  sets $\L_k^{\s_{n-1}}$  and  $\Int \L_k^{\s_{n-1}}$ 
can be divided into the disjoint union of the part in $\L_k$ and the part in the complement of $\L_k$.
We now compute these pieces.

\noindent
{\bf Proposition \BB.5a. (The Parts that Lie in $\L_k$).}
{\sl    For all $k=1, ... , n$ we have:

(i) \ \ \ \ $\L_k^{\s_{n-1}} \cap \L_k \ =\ \L_k,$

(ii)\ \ \ $(\Int \L_k^{\s_{n-1}})\cap \L_k \ =\  ( \Int \L_{k+1})\cap \L_k,$

%(iii) \ \  $(\partial  \L_k^{\s_{n-1}}) \cap \L_k \ =\ \L_k \sim \Int \L_{k+1}.$
}

\noindent
{\bf Proof.}  
 By (\BB.12) we have $\L_k\ss \L_k^{\s_{n-1}}$, and so (i)  follows.

For (ii) note first that by (\BB.12) we have $\Int \L_k^{\s_{n-1}} \ss \Int \L_{k+1}$ and therefore that
$(\Int \L_k^{\s_{n-1}})\cap \L_k \ss ( \Int \L_{k+1})\cap \L_k$.  For the equality in (ii),
assume $A\in  ( \Int \L_{k+1})\cap \L_k = \{r_\ell < 0\} \cap \{r_{\ell+1}\leq 0\}$.
By Fact 1 we have $r_\ell \leq c_\ell \leq r_{\ell+1} \leq 0$.  Now $c_\ell=0$ is impossible,
because this would imply $c_\ell = r_{\ell+1}$, and by Fact 2
this would  imply that $r_\ell = c_\ell = r_{\ell+1} = 0$, which contradicts the assumption
that $r_\ell<0$.  This proves that $r_\ell < c_\ell < r_{\ell+1} \leq 0$,
which implies that $A \in \Int \L_k^{\s_{n-1}} = \{c_\ell < 0\}$.   \qed

%Part (iii) is  a consequence of (i) and (ii).  \qed

Now we compute the $(\sim\L_k)$-part of $\L_k^{\s_{n-1}}$ and
$\Int  \L_k^{\s_{n-1}}$.

\noindent
{\bf Proposition \BB.5b. (The Parts that Lie in $\sim \L_k$).}
{\sl    For all $k=1, ... , n$ we have:

(i)\ \ \ \   $\L_k^{\s_{n-1}}  \cap (\sim \L_k) \ =\ (\L_{k+1} \sim \L_k) \cap \{ \s_n \s_j <0\}$

\noindent
where   $j$ is the largest integer $0\leq j\leq n-1$ with $\s_j(A)\neq 0$
 (with $\s_0(A) \equdef 1$).

(ii) \ \ \ $(\Int \L_k^{\s_{n-1}} ) \cap (\sim \L_k) \ =\ (\Int \L_{k+1} \sim \L_k) \cap \{\s_n\s_{n-1} < 0\}.$

%(iii) \ \ \ $\partial  \L_k^{\s_{n-1}} \cap (\sim \L_k) \ =\ ...$
}

\noindent
{\bf Proof.}
We begin with (i).
Since $\L_k^{\s_{n-1}}$ is contained in $\L_{k+1}$, the $\sim \L_k$-part of $\L_k^{\s_{n-1}}$
is equal to its intersection with
$\L_{k+1} \sim \L_k = \{r_\ell\leq 0<r_{\ell +1}\}$.  Since $r_\ell  \neq r_{\ell+1}$ here,
we have from Fact 2 that
$$
\L_{k+1} \sim \L_k  \ =\   \{r_\ell\leq 0<r_{\ell +1}\} \cap  \{r_\ell <  c_\ell < r_{\ell +1}\}.
\eqno{(\BB.13)}
$$
Intersecting with $ \L_k^{\s_{n-1}}$ we have
$$
E_k \equdef   \L_k^{\s_{n-1}} \cap (\sim \L_k)  =  \L_k^{\s_{n-1}} \cap (\L_{k+1} \sim \L_k)
=  \{r_\ell < c_\ell \leq 0 < r_{\ell+1}\}
\eqno{(\BB.14)}
$$
In particular 0 is not a root of $p(t)$ so that 
$$
p(0) \ =\ \s_n(A) \ \neq\ 0 \qquad{\rm for} \ \ A\in  E_k.
\eqno{(\BB.15)}
$$

Now by Fact 3 we have
$$
E_k \ =\  (\L_{k+1}\sim\L_k) \cap \{A : p_A(t)p_A'(t) <0 \ \ {\rm on} \ \ (c_\ell, r_{\ell+1})\}.
\eqno{(\BB.16)}
$$
$$
\begin{aligned}
&\ \ \ \text
{
Set $\s_0(A) \equdef 1$ and let $j_A$ (or $j$) denote} \\
&\text{ the largest integer $0\leq j\leq n-1$ with $\s_j(A)\neq 0$.
}
\end{aligned}
\eqno{(\BB.17)}
$$
Then we have the Taylor expansion
$$
p_A(t) \ = \ \s_n(A) + t^{n-j} \s_j(A) + t^{n-j+1}  \s_{j-1}(A) + \cdots + t^n.
$$
Consequently, $p_A(t)p_A'(t)<0$ for $0<t$ small $\iff$ $\s_n(A)\s_j(A) <0$.  By (\BB.15) we have that
$\s_n(A)\neq 0$.  This proves Part (i).  

Note: If $c_\ell<0$, then $0\in (c_\ell, r_{\ell+1})$ so that  $p(t)p'(t) <0$ on $(c_\ell, r_{\ell+1})$ is equivalent 
to $\s_n(A) \s_{n-1}(A) = p(0)p'(0)<0$.  The expansion   is only needed for the proof if $c_\ell=0$.

    For Part (ii) note that
since $\Int \L_k^{\s_{n-1}}  \ss \Int \L_{k+1}$, we have
$$
\Int \L_k^{\s_{n-1}}  \cap (\sim\L_k) \ =\ \Int \L_k^{\s_{n-1}} \cap  (\Int \L_{k+1} \sim\L_k),
$$
which equals
$$
\{c_\ell <0\} \cap \{r_\ell < 0 < r_{\ell+1} \}.
$$
By Fact 2, $r_\ell\neq r_{\ell+1}$ implies $r_\ell<c_\ell<r_{\ell+1}$.  Therefore,
$$
(\Int \L_k^{\s_{n-1}}) \cap (\sim \L_k) \ =\ \{r_\ell < c_\ell < 0 < r_{\ell+1}\}
$$
This is a subset of $E_k$ defined above. Since by (\BB.16) we have 
$$
E_k \ =\ (\L_{k+1} \sim\L_k)  \cap \{A : p_A(t)p_A'(t) <0\ \ {\rm on}\ (c_\ell, r_{\ell+1})\}
$$
and since $0\in (c_\ell, r_{\ell+1})$ we have  $p_A(0)p_A'(0) = \s_n(A)\s_{n-1}(A) < 0$,
proving that
$$
(\Int \L_k^{\s_{n-1}}) \cap (\sim \L_k) \ \ss \ (\Int \L_{k+1} \sim \L_k) \ \ss\ 
(\Int \L_{k+1}\sim\L_k) \cap \{\s_n\s_{n-1}<0\}.
$$
Now suppose $A$ is an element of the R.H.S., i.e., assume $r_\ell <0<r_{\ell+1}$,
so by Fact 2, $r_\ell<c_\ell < r_{\ell+2}$ and $p(0)p'(0)<0$.  Then by Fact 3, $0\in (c_\ell, r_{\ell+1})$ proving 
that $c_\ell <0$, i.e., that $A\in \Int \L_k^{\s_{n-1}}$.  This proves Part (ii).

Combining (i) and (ii) we have that the piece of the boundary $\partial  \L_k^{\s_{n-1}}$
contained in $\sim \L_k$ is given by (iii).  \qed

%%%%%%%%%%%%%%%%%%%%%%%%%%%%%%%%%%%%%%%%%%%%%%%%
%%%%%%%%%%%%%%%%%%%%%%%%%%%%%%%%%%%%%%%%%%%%%%%%
%%%%%%%%%%%%%%%%%%%%%%%%%%%%%%%%%%%%%%%%%%%%%%%%
%%%%%%%%%%%%%%%%%%%%%%%%%%%%%%%%%%%%%%%%%%%%%%%%
%%%%%%%%%%%%%%%%%%%%%%%%%%%%%%%%%%%%%%%%%%%%%%%%

\def\elll{\ell}

To complete the proof of Theorem \BB.1$'$ we must prove (2)$'$.

\noindent
{\bf Proof of Theorem \BB.1$'$ (2)$'$.}
First we prove that $\Int \L_k^{\s_{n-1}} \ss \AIthk$.
If $A\in \Int \bL_k^{\s_{n-1}}$, then for $\e>0$ small, $A_\e \equiv A-\e I\in \Int \bL_k^{\s_{n-1}}$.
 Proposition \BB.5 part (ii)  gives two possibilities for $A_\e$.
In both cases we will show that $A\in \AIthk$.

First suppose $A_\e$ is in the $\L_k$ piece of $\Int \bL_k^{\s_{n-1}}$.  
Then $A_{\e\over 2} = A_\e + {\e\over 2} I\in \Int \bL_k$.
By (1)$'$, $A_{\e\over 2} \in \AIth$ for all $\theta \in I_k$.  But $\theta_k$ is the left endpoint of $I_k$.
Hence, picking $\theta \in I_k$, we have that $\theta_k<\theta$ so that
$\AIth \ss \Int \abbf_{\theta_k}$, proving that 
$A_{\e\over 2}$, and hence $A = A_{\e\over 2} +{\e\over 2} I$ is in $\AIthk$.

Secondly, suppose  $A_\e$ is in the $\sim \L_k$-piece of  $\Int \bL_k^{\s_{n-1}}$.  
Then by Proposition \BB.5b(ii)
 $A_\e \in (\Int \bL_{k+1}\sim \bL_k) \cap \{\s_n\s_{n-1}<0\}$.
Since  $A_\e \in\bL_{k+1}\sim \bL_k$, 
Lemma \BB.4 says that $q(A)=k$.  Since $\s_{n}(A_\e)  \s_{n-1}(A_\e) < 0$, we have
by Corollary \BB.3(2), that $f(tA_\e)$ strictly decreases to $(n-2q(A_\e))\pitwo = (n-2k)\pitwo$ 
as $t\nearrow \infty$.  
Thus $tA_\e \in \bbf_{\theta_k}$ for all $t>0$, which, 
by Definition \AA.2  of $\AIthk$, proves that $A\in \AIthk$.

It remains to show that $\AIthk \ss \Int \bL_k^{\s_{n-1}}$, so suppose $A\in \AIthk$.  
Then by definition we have $\liminf_{t\to\infty} f(tA_\e) \geq \theta_k$ for some $\e>0$ small.
By (\BB.7) we can consider $\e>0$ with
 $\s_{n-1}(A_\e) \neq 0$ and $\s_{n}(A_\e) \neq 0$, and apply Lemma \BB.2
to conclude that
$$
\theta_k \ \equiv\ (n-2k)\pitwo \ \leq \ \liminf_{t\to\infty} f(tA_\e) \ =\ (n-2q(A_\e))\pitwo.
$$
This is one of the equivalent ways of saying that $A_\e \in \bL_{k+1}$. 
By decreasing $\e$ we have $A_\e\in\Int \L_{k+1}$.
 If $A_\e\in\bL_k\cap (\Int \L_{k+1})$, then $A\in \Int \bL_k^{\s_{n-1}}$ by Proposition \BB.5a(ii).
   Otherwise, $A_\e \in \Int \bL_{k+1}\sim\bL_k$, so by Lemma \BB.4, $q(A_\e)=k$.
  By Corollary \BB.3(2), if $\s_{n-1}(A_\e)$ and $\s_{n}(A_\e)$ have the same sign,
  then $f(tA_\e)$ strictly increases to $(n-2k) \pitwo = \theta_k$.  This contradicts Definition \AA.2  
  which says that $A\in \AIthk \iff f(tA_\e) \geq \theta_k \ \forall\, t\geq t_0$.
  Therefore, $\s_{n-1}(A_\e)$ and $\s_{n}(A_\e)$ have opposite  signs, which proves that 
  $A\in \Int \bL_k^{\s_{n-1}}$ by Proposition \BB.5b(ii).

  This completes the proof of Theorems \BB.1$'$ and \BB.1. \qed

  From Theorem \BB.1 and Proposition \BB.5 we have the following.

\Theorem {\BB.6} {\sl
The interior of $\bL_k^{\s_{n-1}}$ is given as follows.

(1)\ \ If $\theta\in I_k$ ($k=1, ... ,n$), then
$$
\Int \abbf_\theta\ =\  \Int \bL_k \ \equiv\ \{A\in\Symn : \l_k(A) > 0\}.
$$

 (2) \ \ If $\theta_k$ ($k=1, ... ,n-1$) is a special value, then
$$
\Int \abbf_{\theta_k}\ =\ \{\Int \bL_{k+1} \cap \bL_k\} \cup E_k^\star
$$
where
$$
E_k^\star \ \equiv \ 
 \left(\Int \bL_{k+1} \sim \bL_{k}\right) \cap\left\{  \s_{n-1}(A)\,\s_{n}(A) < 0\right\}
$$
}

\noindent
{\bf Second proof of Proposition \BB.5 part (i)}  We give here a alternative proof  which is based on  Lemma A.1 in [HP2].
Let 
$$
p(t) \ = \ \prod_{k=1}^n (t+\l_k(A))  \ =\ t^n + \s_1(A) t^{n-1} + \s_2(A) t^{n-2} \cdots + \s_n(A)
$$
be the monic polynomial associated to $A\in\Symn$ as above.  Define $\s_0(A) =1$ and set $\s(A) \equiv (1, \s_1(A), ... , \s_n(A))$.
Let $\wt{\s(A)}$ denote the possibly shorter tuple where all the zero entries are removed.  Define
$$
{\rm Var}(\wt{\s(A)})  \ \equiv\ \text{the number of (strict) sign changes in $\wt{\s(A)}$.}
%\eqno{(\BB.8)}
$$

\Lemma{A.1} {\sl
The $k^{\rm th}$ branch of $\{p=0\}$ equals $\{A: {\rm Var}(\wt{\s(A)}) \leq k-1\}$.
In particular, with $p(t) \equiv \det (tI+A)$,}
$$
\L_{k+1}\sim\L_k \ =\ \{A :  {\rm Var}(\wt{\s(A)}) = k\}.
\eqno{(\BB.18)}
$$

This Lemma also applies to the monic polynomial 
$$
{1\over n} p'(t) \ \equiv \ t^{n-1} + {n-1\over n} \s_1(A) t^{n-2}  + \cdots + {1\over n} \s_{n-1}(A)
$$
proving that 
$$
\begin{aligned}
\L_k^{\s_{n-1}} \ &\equiv\ \{A  :  {\rm Var}(\wt{\a(A)}) \leq k-1\} \\
 &{\rm where}\ \ 
\a(A) \equiv (1, \s_1(A), ... , \s_{n-1}(A)).
\end{aligned}
\eqno{(\BB.19)}
$$
(The coefficients in this polynomial ${1\over n} p'(t)$ have been normalized in the definition of $\a(A)$, without changing their signs or their being zero.)
Note that 
$$
\wt{\s(A)} \ =\ 
\left\{ 
\begin{array}{c}
\wt{\a(A)}\quad {\rm if}\ \ \s_n(A)=0
\\
(\wt{\a(A)}, \s_n(A)) \quad {\rm if}\ \ \s_n(A) \neq 0
\end{array}
\right .
\eqno{(\BB.20)}
$$
In particular, ${\rm Var}(\wt{\a(A)}) \leq {\rm Var}(\wt{\s(A)}) \leq {\rm Var}(\wt{\a(A)}) + 1$.
This proves  that $\L_k \ss \L_k^{\s_{n-1}} \ss\L_{k+1}$ and hence Proposition \BB.5a(i).
To prove Proposition \BB.5b(i) that  $\L_k^{\s_{n-1}} \cap (\L_{k+1}\sim\L_k)
= \{A : \s_n(A)\neq 0 \ {\rm and}\  \s_j(A)\s_n(A) < 0\}$ first note that
if $A\in \L_k^{\s_{n-1}}$, i.e.,  ${\rm Var}(\wt{\a(A)}) \leq k-1$, then in order for 
$A\in \L_{k+1} \sim \L_k$, i.e., ${\rm Var}(\wt{\s(A)})  = k$, one must have $\s_n(A) \neq 0$.

Now $\wt{\a(A)} = (1,..., \s_j(A))$ where $j$ is the largest integer $\leq n-1$ such that $\s_j(A)\neq 0$
and $\s(A) = (1, ... , \s_j(A), \s_n(A))$ if $A\in \L_{k+1}\sim\L_k$.
Now $A\in  \L_k^{\s_{n-1}} \cap (\L_{k+1} \sim \L_k) \iff   \s_n(A) \neq 0$ and 
$ {\rm Var}(\wt{\s(A)}) =  {\rm Var}(\wt{\a(A)}) +1$, 
so that $\s_n(A)\neq0$, $\s_j(A)\neq 0$ and $\s_j(A)\s_n(A)<0$.\qed

\vskip .3in
%\vfill\eject

%%%%%%%%%%%%%%%%%%%%%%%%%%%%%%%%%%%%%%%%%%%%%%%%
%%%%%%%%%%%%%%%%%%%%%%%%%%%%%%%%%%%%%%%%%%%%%%%%
%%%%%%%%%%%%%%%%%%%%%%%%%%%%%%%%%%%%%%%%%%%%%%%%
%%%%%%%%%%%%%%%%%%%%%%%%%%%%%%%%%%%%%%%%%%%%%%%%
%%%%%%%%%%%%%%%%%%%%%%%%%%%%%%%%%%%%%%%%%%%%%%%%
%%%%%%%%%%%%%%%%%%%%%%%%%%%%%%%%%%%%%%%%%%%%%%%%
%%%%%%%%%%%%%%%%%%%%%%%%%%%%%%%%%%%%%%%%%%%%%%%%

\noindent {\headfont \CC. The Refined  Dirichlet Problem for $\bbf_\theta$}     

Theorem \AA.12 recalled the general result which applies to the Dirichlet problem
for the equation $\bbf_\theta$.  In this section we will get a deeper result
by using our explicitly computed asymptotic  interiors for $\bbf_\theta$.

We begin with the following note.  In general  for the Dirichlet problem we require
that the boundary be  strictly $\bbf$ and $\wt \bbf$ convex at each point.
(See Theorem \AA.12 for example.)  However, for the subequation $\bbf_\theta$
only one of these is necessary.
Since $\arctan x$ is odd, it is easy to compute that the dual
$$
\wt {\bbf_\theta} \ =\ \bbf_{-\theta},\ \ \ {\rm and\ hence}
\eqno{(\CC.1)}
$$
$$
\oa{{\wt {\bbf_\theta}} } \ =\  \oa {\bbf_{-\theta}}.
\eqno{(\CC.2)}
$$
(Recall the definition of strict $\bbf$-convexity from (\AA.9).)

If $0\leq \theta < n{\pi\over 2}$, then 
$$
\bbf_{-\theta} \supset \bbf_\theta, \ \ \text{and hence}\ \ \oa{{\wt {\bbf_\theta}}}  = 
 \oa {\bbf_{-\theta}} \supset  \oa {\bbf_{\theta}},
\eqno{(\CC.3)}
$$
that is, 
$$
\text
{
$\bo$ is strictly $\bbf_\theta$-convex 
$\qquad\Rightarrow\qquad$
$\bo$ is strictly $\wt\bbf_\theta$-convex 
}
\eqno{(\CC.4)}
$$
Hence, 
$$
\begin{aligned}
&\text{
in applying  Theorem \AA.11 or Theorem \AA.12 to $\bbf=\bbf_\theta$,} \\
&\text
{
one need only hypothesize that $\bo$ is strictly $\bbf_{|\theta|}$-convex.
}
\end{aligned}
\eqno{(\CC.5)}
$$

Note that since $\wt {\bbf_\theta} \ =\ \bbf_{-\theta}$, we may assume in analyzing the 
Dirichlet problem, that $\theta\geq 0$.  For if $\theta<0$, we simply consider the equivalent
problem given by replacing $\bbf_\theta$ with  $\wt {\bbf_\theta}$.

Combining these remarks and Theorem \BB.1 with Theorem \AA.12
yields one of the main results of this paper -- part (b) below.  Part (a) follows from Theorem \AA.10.

\Theorem{\CC.1}  {\sl

(a) Comparison holds for $\bbf_\theta$ for all $\theta\in (-n{\pi\over 2}, n{\pi\over 2})$ on all
bounded domains $\O\ss\rn$.

(b) Suppose $\O$ has a smooth boundary $\bo$, and assume w.l.o.g. that $\theta\geq 0$.

\qquad Case (1): \  If  $\theta\in I_k$ for some integer  $k$,  \\
\ \medskip\hskip 1in then  assume that $\bo$ is strictly 
$
\abbf_{\theta} \ =\ \bL_k$-convex.
 
\qquad Case (2): \ If  $\theta = \theta_k$ is a special value,  \\
\ \medskip\hskip 1in then  assume that $\bo$ is strictly 
$
\abbf_{\theta_k}\ =\ \bL_k^{\s_{n-1}}
$-convex.

Then Perron existence holds for $\bbf_\theta$ on $\O$.
}

We remind the reader of the geometric characterizations of these
boundary convexity hypotheses given in Theorem \AB.1.

\vskip .3in
%\vfill\eject

%%%%%%%%%%%%%%%%%%%%%%%%%%%%%%%%%%%%%%%%%%%%%%%%
%%%%%%%%%%%%%%%%%%%%%%%%%%%%%%%%%%%%%%%%%%%%%%%%
%%%%%%%%%%%%%%%%%%%%%%%%%%%%%%%%%%%%%%%%%%%%%%%%
%%%%%%%%%%%%%%%%%%%%%%%%%%%%%%%%%%%%%%%%%%%%%%%%
%%%%%%%%%%%%%%%%%%%%%%%%%%%%%%%%%%%%%%%%%%%%%%%%
%%%%%%%%%%%%%%%%%%%%%%%%%%%%%%%%%%%%%%%%%%%%%%%%
%%%%%%%%%%%%%%%%%%%%%%%%%%%%%%%%%%%%%%%%%%%%%%%%

\noindent {\headfont \DD.  The Inhomogeneous Dirichlet Problem for the  SL Potential Operator}     

We now consider the Dirichlet problem for  $\tr\{\arctan(D^2_x u)\} = \psi(x)$.

Inhomogeneous equations of this sort were studied in [IDP] where the following was introduced.
An operator $f$, such as $f(A) = \tr\{\arctan \,A\}$, is said to be {\bf tamable} on a subequation 
$\bbf$ if there exists $\chi: f(\bbf) \to \bbr$ strictly increasing, such that $g(A) \equiv \chi(f(A))$ is 
{\bf tame} on $\bbf$ where ``tame'' means that for all $t>0$ there exists $c(t)>0$ such that 
$$
g(A+tI)-g(A) \ \geq\ c(t) \qquad \forall\, A\in\bbf.
\eqno{(\DD.1)}
$$
For tamable operators there are nice results for the Dirichlet problem [IDP].   

The following result was inspired by the paper of Collins, Picard and Wu
 [CPW].

\Theorem{\DD.1}  {\sl
If (and only if) the phase $\theta$ belongs to the highest phase interval 
$I_1 = ((n-2) \pitwo, n\pitwo)$, then the SL potential operator 
$f(A) = \tr\{\arctan \,A\}$  is tamable on $\bbf_\theta$.
More specifically, $\tan\{{1\over n} f(A)\}$ is tame on 
$f^{-1}\{[(n-2){\pi\over 2}+\d, n{\pi\over 2})\} \equiv \bbf_\theta$ where 
$\theta \equiv (n-2){\pi\over 2}+\d$ with $\d>0$.
(The ${1\over n}$ factor is a matter of convenience and not necessary here.)
}

Note that $\aFth = \cp$ if $\theta \in I_1$ by Theorem \BB.1, and hence the appropriate 
boundary convexity is just ordinary  strict convexity.  By Theorem 2.7$'$ in [IDP],  Theorem \DD.1 
gives a different proof of Part A of  the following.

\Theorem{\DD.2} 

{\bf A. (S. Dinew, H.-S. Do and T. D. T\^o  [DDT])}
{\sl
Suppose $\O$ is a bounded domain with smooth strictly convex boundary $\bo$.
Then for any inhomogeneous term $\psi \in C(\ob)$ with values 
$\psi(\ob) \ss  ((n-2) \pitwo, n\pitwo)$ and any boundary function $\vf \in C(\bo)$,
there exists a unique solution $u  \in C(\ob)$ to}
$$
\tr\left\{ \arctan ( D_x^2 u)\right\} \ =\ \psi(x) \quad {\rm on}\ \ \O \quad {\rm with}\ \ u\bigr|_{\bo} =\vf.
\eqno{(\DD.2)}
$$

{\bf B. (M. Cirant and K. Payne [CP, Thm.\ 6.18])}
{\sl
Let $\O\ss\ss\rn$ be a domain and consider an inhomogeneous term $\psi \in C(\ob)$ with values 
in $I_k$, i.e.  
$$
\psi(\ob)\  \ss  \ I_k = \ \left((n-2k) {\pi\over 2}, (n-2(k-1){\pi\over 2}\right)
$$
for some $k$ with $1\leq k  \leq n$.     Then comparison holds for the 
subequation  $\bbf_\psi \ss  \O\times \Symn$ given by}
$$
\bbf_{\psi} \ \equiv \ \{(x,A) :  \tr\{\arctan\ A\}  \geq \psi(x)\}
$$

{\bf C.}  {\sl   Let $\O\ss\ss\rn$ be a domain with smooth boundary which is strictly {$\min\{k, n-k+1\}$-convex}
(by Definition \AB.2).  Let $\psi\in C(\ob)$ be as in Part B.   Then there exists a unique solution $u\in C(\ob)$
to the Dirichlet problem (\DD.2)  for all continuous boundary values $\vf\in C(\bo)$.
}

\noindent
{\bf Proof of Part C.}  This follows from Part B together with Theorem 13.3 in [DDR]  and Theorem \AB.1 above.  \qed

Part C is the best known result on the inhomogeneous Dirichlet problem in the continuous case.

A smooth version of Part A was proved (prior to Theorem \DD.2) in [CPW].

\Theorem{\DD.3. (Collins, Picard and Wu)}  {\sl
Suppose $\O$ is a bounded domain with smooth strictly convex $C^4$ boundary $\bo$.
Consider an inhomogeneous term $\psi \in C^2(\ob)$ with values 
$\psi(\ob) \ss I_1 \equiv ((n-2) \pitwo, n\pitwo)$ and any boundary function $\vf \in C^4(\bo)$.
Suppose there exists a function $\underline u \in C^4(\ob)$ which is a subsolution on $\O$
and satisfies $\underline u = \vf$ on $\bo$.  Then  there exists  a unique $C^{3, \a}(\ob)$
solution $u$ to the Dirichlet problem (\DD.2).
 If all data are $C^\infty$, so is $u$.
 }

This problem has a nice geometric meaning, given in the following proposition.
The result was stated in  [CG; (2.19)], with proof left to the reader.  Since this paper is solely
concerned with the SL Potential equation, we have given a proof in Appendix A 
(See Proposition A.1 and Corollary A.2).

\Prop{\DD.4}  {\sl   
Let $L$ be an  Lagrangian submanifold of $\rn\oplus \rn \cong \rn \oplus i\rn = \bbc^n$.
Assume $L$ is a graph over a domain in $\rn\oplus \{0\}$ and is oriented by $x_1, ... , x_n$.
Set $dz \equiv dz_1\wedge \cdots \wedge dz_n$.  Then $$dz\bigr|_L = e^{i\theta} d{\rm vol}_L,$$
and 
$$
\nabla \theta \ =\ -J H
$$
where $H$ is the mean curvature vector of $L$, and $J$ is the complex structure.
}

\Cor{\DD.5}  {\sl Let $u$ be the solution of the inhomogeneous SL equation (5.2) which is $C^2$ 
on $X^{\rm open}\ss\rn$.
  Let $L$ be the graph of $\nabla u$ in $X\times\rn$. Let $\wt \psi$ be the pull-back of  $\psi$  to $L$. Then
$$
\nabla \wt \psi \ =\ - J H.
$$
Thus, a $C^2$ solution of (\DD.2) is giving a Lagrangian submanifold whose mean curvature vector 
$H$ is
$$
H\ =\ -J \nabla \wt \psi  \qquad{\rm on}\ \ \ L
$$
}

\pf  The function $\theta$ is just $\tr\{\arctan D^2 u\}$ pulled back to $L$.\qed

This result was generalized by Micah Warren to a family of operators and associated
Lagrangian graphs which he found in his thesis.  For this see (6.7) in [Wa].

\noindent
{\bf Open Question:}  It is unclear whether or not some version of Theorem \DD.2 remains true for all possible
inhomogeneous terms $\psi\in C(\ob)$, i.e., $\psi$'s taking values in (a compact subset of)
the open interval $(-n\pitwo, n\pitwo)$.

\noindent
{\bf Proof of Theorem \DD.1}  Set $\d \equiv \theta -(n-2)\pitwo > 0$.  
Adopt the notation $h'(A) \equiv {d\over dt} h(A+tI)\bigr|_{t=0}$.
Set 
$$
g(A) \equiv \chi(f(A)), \qquad {\rm where} \ \  \chi(x) \equiv \tan\left({x\over n}\right).
$$
It suffices to show that for some $c>0$
$$
g'(A) \ \geq\ c \qquad\forall\, A\in \bbf_{\theta}.
\eqno{(\DD.3)}
$$

Let $\l_1 \leq \cdots \leq \l_n$ denote the ordered eigenvalues of $A$.
Then $f(A) = \sum_{j=1}^n \arctan(\l_j)$, and hence
$$
f'(A) \ =\ \sum_{j=1}^n {1\over 1+\l_j^2}.
\eqno{(\DD.4)}
$$
Therefore we have
$$
g'(A) \ =\   \chi'(f(A))  \sum_{j=1}^n {1\over 1+\l_j^2}.
\eqno{(\DD.5)}
$$

\noindent
{\bf Case 1: $A>0$.}
Since $0<\l_1\leq \l_2\leq \cdots \l_n$ we have $f(A) =\sum_j \arctan(\l_j) \geq n \arctan(\l_1) >0$.
Now
$$
\chi'(x) \ =\ {1\over n} \left (1+\tan^2\left({x\over n}\right)\right) \ =\ 
{1\over n} \left (1+\chi^2\left({x}\right)\right).
\eqno{(\DD.6)}
$$
Moreover, $\chi'(x)$ is strictly increasing in $x$ for $x\geq 0$.  Hence,
$$
\chi'(f(A)) \ \geq\ \chi'(n\, \arctan \l_1) \ =\ {1\over n} \left(1+ \tan^2(\arctan \l_1)\right) \ =\ {1\over n}(1+\l_1^2).
$$
Hence, we have
$$
g'(A) \ =\ \chi'(f(A))\sum_{j=1}^n {1\over 1+\l_j^2} \ \geq\ {1\over n} (1+\l_1^2)\sum_{j=1}^n {1\over 1+\l_j^2}
\ \geq\  {1\over n}.
$$

\noindent
{\bf Case 2: $\l_1\leq 0$.}
Suppose we can show that $\theta\in I_1$ implies that
$$
-\tan\left(\pitwo-\d\right) \ \leq \l_1(A).
\eqno{(\DD.7)}
$$
Note that
$$
f'(A) \ =\ \sum_{j=1}^n  {1\over 1+\l_j^2} \ \geq \ {1\over 1+\l_1^2}.
$$
Since $\l_1\leq 0$, (\DD.7) implies that $\l_1^2 \leq \tan^2(\pitwo-\d)$, and hence
$1/(1+\l_1^2) \geq 1/(1+\tan^2(\pitwo-\d))$.  Therefore,
$$
f'(A) \ \geq \ {1\over 1+\tan^2\left(\pitwo-\d\right)}.
\eqno{(\DD.8)}
$$
Note that by (\DD.6) we have $\chi'(x) \geq {1\over n}$.
Therefore, 
$$
g'(A) \ =\ \chi'(f(A)) f'(A) \ \geq\ {{1\over n} \over 1+\tan^2\left(\pitwo-\d\right)}.
$$
In Case 1 the lower bound ${1\over n}$ is larger than here.  Therefore (\DD.3) follows with
$c^{-1}\equiv n(1+\tan^2(\pitwo- \d))$ if we show (\DD.7).  This is immediate from the following.

\Lemma{\DD.6} {\sl
If $A\in \bbf_\theta$ and $(n-2)\pitwo <\theta < n\pitwo$, i.e., $\theta \in I_1$,
then

\noindent
(1)\ \   $A$ must have at least $n-1$ strictly positive eigenvalues, and

\noindent
(2)\ \ if $A$ has a negative eigenvalue, then with $\d\equiv \theta -(n-2)\pitwo$}
$$
-\tan\left(\pitwo-\d\right) \ <\ \l_1(A).
$$
\noindent
{\bf Proof of (1).}  Let $p(A) \equiv \#\{\l_j>0\}$.
Then
$$
(n-2)\pitwo < \theta \leq \sum_{j=1}^n \arctan \, \l_j \leq \sum_{\l_j>0} \arctan\, \l_j < p(A) \pitwo.
$$
Hence, $n-2 <p(A)$.\qed

\noindent
{\bf Proof of (2).}  Note that
$$
(n-2)\pitwo +\d \equiv \theta \leq \sum_{\l_j>0} \arctan \, \l_j - \arctan(-\l_1) < (n-1)\pitwo -\arctan(-\l_1),
$$
and hence $\arctan(-\l_1) < \pitwo-\d$ so that $-\l_1 < \tan(\pitwo-\d)$.\qed

This complete the proof of the ``if'' part of Theorem \DD.1.

\noindent
{\bf Note:}   As an immediate consequence of  Lemma \DD.6, if $\theta$ belongs to the top phase interval
    $I_1$, then each 
$\bbf_\theta$-subharmonic function is quasi-convex, in fact $\tan(\pitwo-\d)$-quasi-convex.

The ``only if'' part of Theorem \DD.1 was proved in [IDP, Prop. 6.19].
We include the proof here for the reader's convenience.

\def\bdf{f}
\def\FTH{\bbf_{\theta}}

\Lemma{\DD.7} {\sl
The  operator $f(A)  \equiv \tr \{\arctan(A)\}$ is  not tamable on $\FTH \equiv \{A : f(A)\geq \theta \}$ for $\theta\leq (n-2){\pi\over 2}$.}

\noindent
{\bf Proof.}  
Consider $A$ with 
$\l_1(A) <<0$ and $\l_k(A)>>0$ for $k>1$.  We can always choose these values
so that $\bdf(A) = (n-2){\pi\over 2}$.  As the absolute value of the eigenvalues becomes
very large the derivative of $\bdf(A)$ goes to zero.  Hence, no matter which smooth function $\chi$
one chooses, the composition $\chi\circ \bdf$ will have 
derivatives going to zero at these points, since $\chi'(\bdf(A))$ will not go to $\infty$
unless $\bdf(A)$ goes to ${n\pi\over 2}$. 
\qed

This completes the proof of Theorem \DD.1.\qed

 \def\leftitem#1{\item{\hbox to\parindent{ #1\hfill}}}

\vskip .3in
%\vfill\eject

%%%%%%%%%%%%%%%%%%%%%%%%%%%%%%%%%%%%%%%%%%%%%%%%
%%%%%%%%%%%%%%%%%%%%%%%%%%%%%%%%%%%%%%%%%%%%%%%%
%%%%%%%%%%%%%%%%%%%%%%%%%%%%%%%%%%%%%%%%%%%%%%%%
%%%%%%%%%%%%%%%%%%%%%%%%%%%%%%%%%%%%%%%%%%%%%%%%
%%%%%%%%%%%%%%%%%%%%%%%%%%%%%%%%%%%%%%%%%%%%%%%%
%%%%%%%%%%%%%%%%%%%%%%%%%%%%%%%%%%%%%%%%%%%%%%%%
%%%%%%%%%%%%%%%%%%%%%%%%%%%%%%%%%%%%%%%%%%%%%%%%

\noindent {\headfont \EE.  A Generalized Version of the Main Theorem \BB.1.}     

The main results in this paper hold with the eigenvalues of $A$ replaced by the G\aa rding eigenvalues of $A$.
We refer to the introduction for definitions and statements.
Let $\ggg$ be  a G\aa rding-Dirichlet polynomial   of degree $m$
 with ordered G\aa rding eigenvalues $\l^\ggg_k$ and associated branches 
 $\bL_k^\ggg \equiv \{ \l_k^\ggg\geq 0\}$, for $k=1,...,m$.
 One has
 $$
\ggg(tI+A) \ =\ \ggg(I) \prod_{k=1}^m (t+\l^\ggg_k(A))
\ =\ \ggg(I) \left( t^m + \s_1^\ggg(A) t^{m-1} + \cdots + \s_m^\ggg(A)      \right).
$$
The derivative
$
\ggg'(A) \equiv {1\over m} {d\over dt} \ggg(tI +A) \bigr|_{t=0}.
$
is  a G\aa rding-Dirichlet polynomial   of degree $m-1$ whose eigenvalues  $\l_k^{\ggg'}$ are the critical 
points of $\ggg(tI+A)$.  Note that
 $$
\s_m(A) \ =\ \ggg(A) \qquad{\rm and}\qquad \s^\ggg_{m-1}(A) \ =\ \ggg'(A).
\eqno{(\EE.1)}
$$
 
 We begin by looking at the asymptotic expansion of the $\ggg$-SL potential operator.
 
 \Lemma{\EE.1}  {\sl
 Suppose $\ggg(A)\neq 0$.  Then 
 $$
 \begin{aligned}
 f^\ggg(tA) \ &\equiv\ \sum_{k=1}^m \arctan \l_k^\ggg(tA) \\
  \ &=\ 
 (m-2q(A))\pitwo \ -  \  {1\over t} {\ggg'(A)\over \ggg(A)} \  +  \ O\left({1\over t^3}\right) \qquad {\rm as}\ \ t\to \infty.
 \end{aligned}
 $$
 where $q(A)$ now  denotes the number of strictly negative $\ggg$-eigenvalues of $A$.}
 
\pf
Note from the first display above that $\l^\ggg_k(tA) = t\l^\ggg_k(A)$ and so 
$ f^\ggg(tA) = \sum_{k=1}^m \arctan ( t \, \l_k^\ggg(A))$. Since $\ggg(A) = \l^\ggg_1(A) \cdots  \l^\ggg_m(A) \neq 0$,
all the eigenvalues are non-zero, and so the difference between the number of strictly positive eigenvalues
and the number of strictly negative ones is $(m-q(A)) -q(A) = m-2q(A)$.  The proof is now the same as 
that of Lemma \BB.2. \qed

\Cor{\EE.2} {\sl
Suppose $\ggg(A)\neq 0$.  Then Corollary \BB.3 holds with $f(tA)$ replaced by $f^\ggg(tA)$.
}

Now Generalized Theorem \BB.1 can be rephrased as Generalized  Theorem \BB.1$'$ exactly as 
Theorem  \BB.1 was rephrased as Theorem  \BB.1$'$.  Also the analogue
of  Lemma \BB.4 holds here. 

The next step is to generalize Proposition \BB.5. Here the statement and arguments are
essentially identical to those  in Section \BB.  We state the part concerning the interior, since that
is what is important for the Dirichlet problem

\Prop{\EE.3}  {\sl
For $k=1, ... , m$, the set $\Int \bL_k^{\s^\ggg_{m-1}}$ is a disjoint union 
$\Int \bL_k^{\s^\ggg_{m-1}}  =(\Int  \bL_{k+1}^\ggg  \cap  \bL_k^\ggg) \cup E_k$ where       
$$
E_k \ \equiv \ (\Int \bL_{k+1}^\ggg  \sim  \bL_k^\ggg ) \cap \{  \s_{m-1}^\ggg(A) \,   \s_{m}^\ggg(A) \ <  0\}.
$$
}

%\pf  One uses Lemma \BB.6 exactly as it is used to prove Proposition \BB.5.\qed

The remainder of the proof of  Generalized Theorem \BB.1$'$  follows exactly the argument given for  the ``Proof of Theorem \BB.1$'$ (2)$'$''
in Section \BB. \qed

\vskip .3in

\centerline{\bf The SL Curvature Operator}  
\centerline{\bf  for the Second Fundamental Form of the Graph}
\smallskip

Here we include some brief remarks on another way to diversify the SL operator
by looking at the second fundamental form of the graph of a scalar function $u$.
For this operator very little is known concerning uniqueness (for the (DP)).
However,  our result establishing the appropriate strict boundary convexity 
carries over and provides existence for the broad class of appropriately ``pseudo-convex'' domains.

\def\slc{SL^{\rm curv}}
\def\Kt{{\mathcal F}_\theta}

\Def{\EE.4.  (The SL-Curvature Operator)}
Given a smooth function $u$ on an open subset $X\ss\rn$, let 
$\kappa_1(x), ... , \kappa_n(x)$ denote the principal curvatures of it graph $M\ss \bbr^{n+1}$.
Replacing the eigenvalues of $D_x^2 u$ in (\II.1), by these principal curvatures, yields the 
{\bf SL-curvature operator}  $\slc$ defined by (\EE.2) below.

This operator was first studied by Graham Smith.  In [GS1] he restricts to functions $u$ with
$D^2 u > 0$  and  $\theta \in [(n-1){\pi\over 2}, n{\pi\over 2})$, 
and he looks at the number $\rho$ such that $\sum_j \arctan\{\rho \kappa_j\}=\theta$.
 When $\theta =  (n-1){\pi\over 2}$, he gives a very nice geometric interpretation
of this ``curvature'' $\rho$ (see (i)--(iii) on page 59 of [GS1]).

As noted in [DDR, Sec. 11.5], using jet variables $p=D_x u$ and $A=D^2_x u$, if we define
$$
E_p \ =\ P_{p^\perp} +{1\over \sqrt{1+ |p|^2}} P_p
\ =\ 
I -  {|p|^2  \over  \sqrt{1+ |p|^2} (1+ \sqrt{1+ |p|^2}) } P_p,
\eqno{(\EE.1)}
$$
then 
$$
II(u) \ \equiv \ {1\over \sqrt{1+ |p|^2}} E_p A E_p
$$
 is the second fundamental form of the graph $M$
of $u$, so that its eigenvalues  are the principal curvatures $\kappa_1, ... , \kappa_n$ of $M$.  Thus
$$
\slc(p, A) \ \equiv\ \tr\left\{ \arctan {II(u) } \right\} \ =\ \sum_{j=1}^n \arctan \kappa_j.
\eqno{(\EE.2)}
$$
For phases $\theta$ in the allowable range $- n {\pi\over 2}  < \theta < n {\pi\over 2}$, 
$$
\Kt \ \equiv\ \{(p,A) : \slc(p,A) \geq \theta \}
\eqno{(\EE.3)}
$$
is a subequation, which is constant coefficient and reduced, but {\bf not} pure second order.
The positivity requirement,  that $\slc(A+P) \geq \slc(A)$ for $P\geq 0$,  follows from the fact that 
$\bra {E_p P E_p x} x = \bra {P(E_p x)} {E_px}$ so that $P\geq 0 \Rightarrow E_p P E_p \geq 0$.

In fact for each fixed  $p\in\rn$,
$$
g(A) \ \equiv\ \det\left( {1\over \sqrt{1+ |p|^2}} E_p A E_p    \right) \ =\ \det(II(u))
\eqno{(\EE.4)}
$$
is a G\aa rding/Dirichlet polynomial in $A$ whose eigenvalues 
 (the negatives of the roots of $\det\left(tI +E_p A E_p    \right)$) 
 are the  curvatures $\kappa_1(p,A), ... , \kappa_n(p, A)$.
In particular, $\Kt$ fibres over $\rn$, with fibre at $p$ given by  the pure second-order subequation
$$
\begin{aligned}
(\Kt)_p \ &=\   {1\over {1+ |p|^2}}    E_p \bbf_\theta E_p, \ \ \ {\rm i.e.} \\
 (p,A) \in \Kt \ \ &\iff\ \ B\equiv (1+|p|^2) E_p^{-1} A E_p^{-1} \in \bbf_\theta.
\end{aligned}
\eqno{(\EE.5)}
$$
In addition, the fibres of the interior are the interiors of the fibres.  Consequently,
the asymptotic interior $\oa \Int \Kt$ can be computed fibrewise with the answer given by Theorem \BB.1$'$.
This gives the following optimal pseudoconvexity result for the SL curvature operator.
The  proof is left to the reader.

\noindent
{\bf   Theorem \BB.1$'$ for the SL Curvature Operator.}      
Let ${\mathcal L}_k$ be the set where the $k^{\rm th}$ ordered eigenvalue $\kappa_k(p,A)$
  is $\geq 0$.  Let  ${\mathcal L}_k^{\s_{n-1}}$ also be defined in analogy with section 3.

(1)$'$ \ \ If $\theta\in I_k$ ($k=1, ... ,n$), then
$$
 \oa{\Int}\Kt\ =\  \Int {\mathcal L}_k.
$$

(2)$'$ \ \ If $\theta_k$ ($k=1, ... ,n-1$) is a special value, then
$$
\oa{\Int} \, {\mathcal F}_{\theta_k}\ =\  \Int {\mathcal L}_k^{\s_{n-1}}.
$$

Given $\theta\in (- n {\pi\over 2}, n {\pi\over 2})$, consider the (DP) for $\Kt$ on a domain
$\O\ss\ss \rn$ with general boundary function $\vf\in C(\bo)$.  None of the available comparison
techniques seem to apply.  Thus:
$$
\text{Comparison for $\Kt$ remains an interesting open question.}
\eqno{(\EE.6)}
$$

What we do know from [DDR] can be outlined as follows.  From 
 Definition 8.2 and Theorem 10.1 in [DDR] we have:
$$
\text{A weak form of comparison holds.}
\eqno{(\EE.7)}
$$

Assume that $\bo$ is smooth and strictly ${\mathcal F}_{|\theta|}$-convex.  Then from 
Theorem 13.4 in [DDR]:
$$
\text{Existence holds for the $\Kt$ (DP) on $\O$ for all $\vf\in C(\bo)$,}
\eqno{(\EE.8)}
$$
and from Theorem 12.7 in [DDR]:
$$
\begin{aligned}
&\qquad \text{All solutions for the $\Kt$  (DP) on $\O$ for $\vf\in C(\bo)$.}  \\
&\text{are squeezed between a maximal solution and a minimal solution,}  \\
&\qquad\qquad \ \ \  \text{  namely the two Perron solutions.}
\end{aligned}
\eqno{(\EE.9)}
$$

\Remark{\EE.5. (Radial Harmonics)}
Suppose $u(x)\equiv \psi(|x|)$, $r\equiv |x|$ is a radial solution to $\slc(D_x u, D_x^2 u) =\theta$.
Then $p=D_x u = \psi'(r) {x\over r}$ and $A = D_x^2 u = {\psi'(r) \over r} P_{x^\perp} + \psi''(r) P_x$.
Note that $A$ commutes with $E_p$.  Set $y(r) \equiv \psi'(r)$.  One computes that 
$$
{1\over \sqrt{1+|p|^2} } E_p A E_p \ =\ 
{1\over  \sqrt{1+y^2} } {y\over r} P_{x^\perp}
 + {y'\over (1+y^2)^{{3\over 2}} }P_x.
\eqno{(\EE.10)}
$$
Hence the principal curvatures of the graph of $\psi(|x|)$ are: 
$$
{1\over r} {y\over   \sqrt{1+y^2} } \ \ \text{ with multiplicity $n-1$ and }
 {y'\over (1+y^2)^{{3\over 2}} } \ \ \text{ with multiplicity 1}.
 $$
It follows that:
$$
\begin{aligned}
0\ &=\ {\rm Im}\left(e^{-i\theta} \left( 1 + i {y \over r\sqrt{1+y^2} }         \right)^{n-1}
  \left( 1 + i {y' \over (1+y^2)^{{3\over 2}}   }      \right)        \right)    \\
  & =   \ {1\over nr^{n-1}}  {d\over dr} {\rm Im} \left(e^{-i\theta} \left( r + i {y \over \sqrt{1+y^2} }         \right)^{n}\right)
  \end{aligned}
\eqno{(\EE.11)}
$$
since $ {d\over dr} {y \over \sqrt{1+y^2} }  = {y' \over (1+y^2)^{{3\over 2}}} $.

This proves that $\slc(D_x u, D_x^2 u) = \theta$ implies
$$
{\rm Im} \left(e^{-i\theta}  \left( r + i {y \over \sqrt{1+y^2} }         \right)^{n}\right) \ =\ c \qquad {\rm a \  constant}.
\eqno{(\EE.12)}
$$
Compare this with the calculation in [CG] pp. 99-100 that if $SL(D_x^2 u) = \theta$, then Im$(r+iy)^n = c$, a constant.

We leave it to the reader to calculate some low dimensional cases.

\vskip .3in
%\vfill\eject

%%%%%%%%%%%%%%%%%%%%%%%%%%%%%%%%%%%%%%%%%%%%%%%%
%%%%%%%%%%%%%%%%%%%%%%%%%%%%%%%%%%%%%%%%%%%%%%%%
%%%%%%%%%%%%%%%%%%%%%%%%%%%%%%%%%%%%%%%%%%%%%%%%
%%%%%%%%%%%%%%%%%%%%%%%%%%%%%%%%%%%%%%%%%%%%%%%%
%%%%%%%%%%%%%%%%%%%%%%%%%%%%%%%%%%%%%%%%%%%%%%%%
%%%%%%%%%%%%%%%%%%%%%%%%%%%%%%%%%%%%%%%%%%%%%%%%
%%%%%%%%%%%%%%%%%%%%%%%%%%%%%%%%%%%%%%%%%%%%%%%%

\noindent {\headfont \FF.  Results on Riemannian Manifolds.}     

Via the work in [DDR] the results above can be carried over to fairly general spaces.
Let $\O\ss\ss X$ be a domain with smooth boundary in a riemannian $n$-manifold, 
and let $\Hess \, u \in C^\infty(\Sym(T^*(X))$ be the riemannian Hessian.  Then for
$C^2$-functions $u$ the SL potential operator
$$
f(\Hess \, u) \ =\ \tr\{\arctan(\Hess \,  u)\}
$$
makes sense and extends to upper semi-continuous functions.  For $\theta\in (-n\pitwo, n\pitwo)$
we have the subequation $\bbf_\theta$ on $X$ and the associated equation $\partial \bbf_\theta$.
These equations are locally jet-equivalent to the constant coefficient equations discussed above,
and so the work in [DDR] applies.  

The strict  boundary convexity assumption on the second fundamental
form, analyzed in Section \BB, carries over directly to $\bo$.  We assume this is satisfied, and that
 there exists a smooth convex function defined on $X$.

\Theorem{\FF.1. [DDR]}  {\sl
For each
$\vf\in C(\bo)$ there exists a unique solution $u\in C(\ob)$ to the Dirichlet problem, i.e., 
$u$ is an $\bbf_\theta$-harmonic function on $\O$ and $u\bigr|_{\bo} = \vf$.
}

This extends more generally as follows.  Let $X$ be a riemannian $n$-manifold with a topological
$G$-structure for a compact group $G\ss {\rm O}(n)$.  Let $\ggg$ be a $G$-invariant G\aa rding-Dirichlet polynomial
of degree $m$ on $\Symn$.  Then we have a well defined SL potential operator
$$
f^\ggg(\Hess \, u) \ =\  \sum_{k=1}^m \arctan \{ \l_k^\ggg(\Hess \,  u)\}
$$
where $\l^\ggg_k(A)$ are the G\aa rding eigenvalues of $A$.  We have the subequation $\bbf^\ggg_\theta$ for 
$\theta \in (-m\pitwo, m\pitwo)$ and its associated equation.  We suppose that $X$ carries a smooth strictly
$\G$-subharmonic function where $\G$ is the G\aa rding cone for $\ggg$.

\Theorem{\FF.2}  {\sl
Suppose that $\bo$ satisfies the strict boundary convexity hypothesis for $\bbf_\theta^\ggg$.  Then for each
$\vf\in C(\bo)$ there exists a unique solution $u\in C(\ob)$ to the Dirichlet problem, i.e., 
$u$ is an $\bbf_\theta^\ggg$-harmonic function on $\O$ and $u\bigr|_{\bo} = \vf$.
}

\Ex{\FF.3}  Suppose $(X,J)$ is an almost complex riemannnian manifold where $J$ is orthogonal.
Then we can take $\ggg(A) = \det_{\bbc}(A_\bbc)$ where $A_\bbc \equiv \half (A-JAJ)$.
Here the $\G$-subharmonic functions are exactly the plurisubharmonic functions.
This gives solutions to the complex SL potential equation.  There is a quaternionic analogue.
One also has in the complex case the Lagrangian Monge-Amp\`ere operator discussed in [Lag].

%\vfill\eject
\vskip.3in

%%%%%%%%%%%%%%%%%%%%%%%%%%%%%%%%%%%%%%%%%%%%%%%%%%%%
%%%%%%%%%%%%%%%%%%%%%%%%%%%%%%%%%%%%%%%%%%%%%%%%%%%%
%%%%%%%%%%%%%%%%%%%%%%%%%%%%%%%%%%%%%%%%%%%%%%%%%%%%
%%%%%%%%%%%%%%%%%%%%%%%%%%%%%%%%%%%%%%%%%%%%%%%%%%%%
%%%%%%%%%%%%%%%%%%%%%%%%%%%%%%%%%%%%%%%%%%%%%%%%%%%%
%%%%%%%%%%%%%%%%%%%%%%%%%%%%%%%%%%%%%%%%%%%%%%%%%%%%
%%%%%%%%%%%%%%%%%%%%%%%%%%%%%%%%%%%%%%%%%%%%%%%%%%%%

\centerline{\bf Appendix A.  A Geometric Interpretation of the Inhomogeneous DP.}

The   equation (A.1) below  appeared as equation (2.18) in [CG]. 
 We left the proof as a exercise for the reader.
  However,  this equation has an immediate consequence for the  
 Dirichlet problem for the inhomogeneous SL equation (A.2), which is discussed in Section 5.
This is given in Corollary A.2.  It may have gone unnoticed and   seems not to be well understood.
For the convenience of the reader we give the proof of equation (2.18) in [CG] here.

\Prop {A.1} {\sl
  Let $X$ be a Calabi-Yau manifold of complex dimension $n$.  Let $\Phi$ be the
  parallel $(n,0)$-form normalized so that ${\rm Re} \Phi$ has comass 1.  
Given  $L\ss X$  an oriented Lagrangian
  submanifold, define the phase $\theta$ mod $2\pi$ by  
$$
\Phi \bigr|_L  \ =\  e^{i\theta} d{\rm vol}_L.
$$  
Then 
for any tangent vector field $V$ on $L$, we have
$$
V \theta \ =\ \bra {JV} H,
$$
that is
$$
\nabla \theta \ =\ -JH
\eqno{(A.1)}
$$
where $H$ is the mean-curvature vector field of $L$, and $J$ is the complex structure on $X$.}

Proposition A.1 has the following immediate implication for the inhomogeneous SL potential equation 
$\tr \left \{ \arctan(D^2_x u)  \right \} \ =\ \psi(x)$.  Let  $z\equiv x+iy \in \rn\oplus i \rn = \bbc^n$ 

\Cor {A.2}  {\sl Suppose $L\equiv \{(x, \nabla u(x)) : x\in \O\}$ is the graph
of the gradient of $u\in C^2(\O)$ over a domain $\O\ss \rn$.  Then the inhomogeneous
term
$$
\theta(x) \ \equiv \ \tr\left\{ \arctan D_x^2 u    \right\},
\eqno{(A.2)}
$$
considered as a function on $L$, is the phase function for $L$.
Thus it has gradient related to the mean curvature vector field $H$ of $L$ by
$$
\nabla \theta \ =\ -JH \qquad {\rm on}\ \ L.
\eqno{(A.3)}
$$
Otherwise said, if $u$ is  a solution to the equation
$$
 \tr \left \{ \arctan(D^2_x u)  \right \} \ =\ \psi(x)
$$
on $\O$, with $\psi(x)$ smooth, then 
$$
\nabla \wt\psi \ =\ -J H  \qquad {\rm on}\ \ \  L
\eqno{(A.3)}
$$
where $\wt\psi$ is the pull-back of $\psi$ to $L$.}

\noindent
{\bf Note.}  Proposition A.1 is actually  independent of the orientation of $L$.  A change of orientation
changes the function $\theta$ to $\theta +\pi$, and the conclusion is the same.   In Corollary A.2, $L$ is given the orientation of $\O$.

\noindent
{\bf Proof of Proposition A.1.}  By a complex  linear change of coordinates we may assume at $p$ we have 
  $\Phi = dz_1\wedge \cdots \wedge dz_n$.  Now   let $p=(x_0, \nabla u(x_0))$. 
The map $D^2 u$ is symmetric, so by a change of variables $(x,y) \to (gx, gy)$ for some $g\in {\rm SO}(n)$,
we can assume that at $x_0$, $D^2 u$ is diagonal, i.e.,  $(D^2_{x_0} u)(\e_k) = \l_k \e_k$ for an orthonormal basis
$\e_1, ... , \e_n$ of $\rn$.

Now let $e_1, ... , e_n$ be an oriented orthonormal frame field
 on $L$ in a neighborhood of $p=(x_0, \nabla u(x_0))$.
Then
$$
\Phi (e_1\wedge \cdots \wedge e_n)  
\ =\  e^{i\theta}, 
$$
and so 
$$
V e^{i\theta} \ =\ e^{i\theta} i V\theta \ 
=\ \Phi \left( \sum_{k=1}^n e_1\wedge \cdots \wedge (\nabla_V e_k) \wedge \cdots \wedge e_n  \right )
$$
since $\Phi$ is parallel.
We may assume  $\nabla^L_{e_j} e_k = (\nabla_{e_j} e_k)^{{\rm tang}} = 0$ at $x_0$,  so $\nabla_V e_k =  ( \nabla_V e_k)^{\rm normal}
= B_{V, e_k}$ the second fundamental form of $L$ at $x_0$.  Therefore  we have
$$
V e^{i\theta} \ =\ e^{i\theta} i V\theta \ 
=\  \Phi \left( \sum_{k=1}^n e_1\wedge \cdots \wedge B_{V, e_k} \wedge \cdots \wedge e_n  \right ),
$$
and we can write
$$
B_{V, e_k} \ =\ \sum_{j=1}^n \bra {B_{V, e_k} } {Je_j} Je_j.
$$
Now pick the frame field   at $x_0$ to be 
  $$
  e_k =   {1\over \sqrt{1+\l_k^2}} (\e_k +\l_k J\e_k), \qquad {\rm for} \ \ k=1, ... , n,
  $$
so that at $x_0$ the vectors $e_k$ and $Je_k$ lie in the k$^{\rm th}$ complex coordinate line.
Recall  that at the point $p$,  $\Phi = dz_1\wedge \cdots \wedge dz_n$. 
% (A little point: the orientation of 
%the $dz_k$'s and the $e_k$'s may not agree. Changing the orientation of $L$ changes
%$\theta$ by $\pi$ and the conclusion is the same.)   
Hence, at $p$
$$
\begin{aligned}
V e^{i\theta} \ &=\ e^{i\theta} i V\theta 
=\  \{ dz\}\left( \sum_{k=1}^n e_1\wedge \cdots \wedge B_{V, e_k} \wedge \cdots \wedge e_n  \right )  \\
&= \ \sum_{k} dz_1(e_1) \cdots dz_k\left( \sum_j \bra {B_{V, e_k} } {Je_j} Je_j \right) \cdots dz_n(e_n)  \\
&= \  \sum_{k} dz_1(e_1) \cdots dz_k\left(\bra {B_{V, e_k} } {Je_k} Je_k \right) \cdots dz_n(e_n)  \\
&\equiv \  \sum_{k} dz_1(e_1) \cdots dz_k\left(\a_k Je_k \right) \cdots dz_n(e_n) \qquad {\rm with}\ \ \a_k \equiv 
 \bra {B_{V, e_k} } {Je_j}  \\
&= \   \sum_{k} dz_1(e_1) \cdots i \a_k dz_k\left(e_k \right) \cdots dz_n(e_n)  \\
&= \ i \sum_k \a_k dz(e_1\wedge \cdots \wedge e_n) \ = \ i \left(\sum_k \a_k\right) e^{i\theta}.
\end{aligned}
$$
Hence, with summation convention, 
$$
\begin{aligned}
V \theta \ &=\  \sum_k \a_k \ =\  \bra {B_{V, e_k} } {Je_k}  \ =\ \bra{\nabla_{e_k} V}  {Je_k} \\
&=\ 
-\bra  V {\nabla_{e_k} Je_k}  \ =\ - \bra V {J\nabla_{e_k} e_k} \ =\ \bra {JV} {B_{e_k, e_k}} \\
& = \ \bra {JV} H.
\end{aligned}
$$
\qed

\vskip.3in
%\vfill\eject

%%%%%%%%%%%%%%%%%%%%%%%%%%%%%%%%%%%%%%%%%%%%%%%%%%%%
%%%%%%%%%%%%%%%%%%%%%%%%%%%%%%%%%%%%%%%%%%%%%%%%%%%%
%%%%%%%%%%%%%%%%%%%%%%%%%%%%%%%%%%%%%%%%%%%%%%%%%%%%
%%%%%%%%%%%%%%%%%%%%%%%%%%%%%%%%%%%%%%%%%%%%%%%%%%%%
%%%%%%%%%%%%%%%%%%%%%%%%%%%%%%%%%%%%%%%%%%%%%%%%%%%%
%%%%%%%%%%%%%%%%%%%%%%%%%%%%%%%%%%%%%%%%%%%%%%%%%%%%
%%%%%%%%%%%%%%%%%%%%%%%%%%%%%%%%%%%%%%%%%%%%%%%%%%%%

\centerline{\bf Appendix B.  Remarks Concerning Convexity}

Part of the point of this appendix is to show that strict $\oa \bbf_\theta$-convexity 
of the boundary $\bo$ is the  right, i.e., borderline condition for the Dirichlet problem.
This goes back to work in [CNS] for  $\theta$ in the highest interval.  We also discuss 
how this convexity relates to convexity of the domain $\O$.

\def\hs{h^\star}
\noindent
Let 

%\noindent
  $\O\ss\ss\rn$ be a domain with smooth boundary $\bo$,

%\noindent
$\bbf\ss\Symn$ be a   subequation,

%\noindent
$\l>0$ be such that $(-\cp) \cap (\l  I+ \bbf) = \emptyset$ (this always exists), and

%\noindent 
$\vf = -{\l\over 2} \|x\|^2\bigr|_{\bo}$.

%\Def{2}  The boundary $\bo$ is strictly $\bbf$-convex if its second fundamental form $B$
%w.r.t.\ the interior normal satisfies $B=A\bigr|_{T_x(\O)}$ for some $A\in \Int \

\Def{B.1} $\O$ is said to be $\bbf(\O)$-{\bf convex} if for every $K\ss\ss\O$,
one has ${\wh K}^{\bbf(\O)}\ss\ss\O$ where 
${\wh K}^{\bbf(\O)}\equdef \{x\in \O : f(x) \leq \sup_K f, \ \forall \, f\in \bbf(\O)\}$

\noindent
{\bf Theorem B.2.}  {\sl
Assume there exists a function $h\in C^2(\ob)$ such that 

(i) \ \ $h\bigr|_{\O} \in \bbf(\O)$,

(ii) \ \ $h\bigr|_{\bo} \ =\ \vf.$

\noindent
Then the domain $\O$ is $\bbf$-convex.   

\def\bbg{{\mathbb G}}

Furthermore, let $\hs \equdef h + {\l\over 2} \|x\|^2$.  Then for every $x\in\bo$ where $(D\hs)_x \neq 0$,
 there exits $A\in ( \l I +  \bbf)$ such that  
$$
A\bigr|_{T_x(\bo)} \ =\ c_x B_x
\eqno{(B.1)}
$$
where $c_x>0$ and $B_x$ is the second fundamental form of $\bo$ at $x$ with respect to the interior normal.

In particular, if $\bbf = \oa \bbg$  for a subequation $\bbg$, then $\bo$ is strictly $\bbg$-convex.
}

\noindent
{\bf Proof.} We begin by proving the second assertion.
Fix $p\in\bo$ and w.l.o.g.\ assume $p$ is the origin. 
Choose coordinates $x= (x', x_n)$ such that in a neighborhood of 0
$$
\O \ =\ \{(x', x_n) : x_n \geq g(x')\}
$$
where $g$ is $C^\infty$ with 
$$
g(0) \ =\ 0, \qquad (Dg)_0 = 0, \qquad {\rm and}\qquad  (D^2 g)_0 \ =\ B 
$$
where $B \ \equdef\ \text{the second fundamental form of $\bo$ at 0 w.r.t.\ the interior normal}$. 

In a neighborhood $U$ of  we have a defining  function for $\O$ given by
$$
\rho(x) \ \equiv\ g(x')- x_n
\eqno{(1)}
$$
with
$$
(D^2 \rho)_0 \ \equiv \ B.
$$

\noindent
{\bf Lemma {B.3}.}  {\sl
If  $\wt \rho$ is any other defining function for $\O$ in $U$, then
}
$$
(D^2\wt \rho)\bigr|_{T_0(\bo)} \ \equiv \ |D\wt\rho|_0 B 
$$
\noindent
{\bf Proof.} 
In a  neighborhood of 0 we have that $\wt \rho (x) = a(x)\rho(x)$ where $a>0$.
Now
$$
D\wt\rho \ =\ (Da)\rho + a (D\rho)
$$
$$
D^2\wt\rho \ =\ (D^2 a)\rho + (Da)\circ (D\rho) + a(D^2 \rho).
$$
At $x=0$ we have $\rho(0)=0$ and $(D\rho)_0=(0, ... , 0, -1) \equiv  n$. Therefore,
$$
(D\wt\rho)_0 \ = \ (0, ... , 0 ,  -a(0))
$$
$$
(D^2\wt\rho)_0 \ =\ (Da \circ n) + a(0)(D^2 \rho)_0 \qquad\text{and so}
$$
$$
(D^2\wt\rho)\bigr|_{T_0(\bo)}  \ =\ a(0)(D^2 \rho)\bigr|_{T_0(\bo)}  \ =\ a(0) B.  \hskip .5in\mathqed
$$

Now $\vf = -{\l\over 2}\|x\||^2\bigr|_{\bo}$, and  $ h\in C^2(\ob)$ is $\bbf$-subharmonic
 on $\O$  with boundary values $\vf$.  That is, 

(1) \ \ $D^2 h \ \in \ \bbf$\ \  on \ \ $\ob$,

(2) \ \ $( h-\vf)\bigr|_{\bo}\ =\ 0.$

\noindent
Recall that 
$
h^\star \ \equiv \ h   +  \frac \l 2 \|x\|^2
$
Note that  $h^\star \in C^2(\ob)$ is an $\bbf^\star$-subharmonic function for the subequation
$$
\bbf^\star \ \equdef \ \l I+\bbf
$$
with
$$
\hs\bigr|_{\bo} \ = \ 0.
$$
Note also that 
$
D^2 \hs \ \in\   \l  I+\bbf \ \ss\  \Int \bbf
$
by the positivity condition. Hence, $\hs$ is strictly $\bbf$-subharmonic.
%Furthermore, $\hs$ is $C^2$ up to the boundary $\bo$.
%(By Whitney extension, $\hs$ extends to a $C^2$ function in a neighborhood 
%of $\ob$.)

With the supposition that
$
(D\hs )_x   \neq 0
$.
we have that  $\hs$ is a  defining function for
$\bo$ in a neighborhood of $x$.  We now apply Lemma B.3 to establish the second assertion.
%Applying this to an asymptotic cone $\oa\bbf$ gives the third assertion.

 Now for the first assertion.
 By the definition of $\l$ and Theorem 3.1 in [SMP] we know that $\hs$ satisfies
the Strong Maximum Principle, and so $\hs<0$ on $\O$.  
Suppose  $K\ss\ss\O$ and let $\d \equiv {\rm dist}(K,\bo) <0$.  Then we have ${\wh K}^{\bbf(\O)} \ss \O_\d
\equdef \{x\in \O: \hs(x)\leq \d\}$.  \qed

\vskip .3in
%\vfill \eject

\centerline{\bf References}

\noindent
[BW] S. Brendle  and M.  Warren,  {\sl A boundary value problem for minimal Lagrangian
graphs}, Journal of Differential Geometry {\bf 84} (2010), 267-287.

\noindent
{[CNS]}   L. Caffarelli, L. Nirenberg and J. Spruck,  {\sl
The Dirichlet problem for nonlinear second order elliptic equations. I: Monge-Amp\`ere equation},  
Comm. Pure Appl. Math.  {\bf 37} (1984),   369-402.

\noindent
{[GC1]}  G. Chen,  {\sl On J-equation}, ArXiv preprints (2019), ArXiv:1905.10222.

\noindent
{[GC2]}    ------------- ,   {\sl  Supercritical deformed Hemitian-Yang-Mills equation}, ArXiv preprints (2020), ArXiv:2005.12202.

\noindent
{[CSY]}  J. Chen, R. Shankar and Y. Yuan, {\sl Regularity for convex viscosity solutions
of the special Lagrangian equation}, ArXiv:1911.05452.

\noindent
[CDS1]  X.-X. Chen, S. Donaldson, S. Sun, {\sl  K¬ahler-Einstein metrics on Fano manifolds. I: 
Approximation of metrics with cone singularities}, J. Amer. Math. Soc. {\bf 28} (2015), no.
1, 183Ð197.

\noindent
[CDS2]    ------------- ,  {\sl  K¬ahler-Einstein metrics on Fano manifolds. II:
Limits with cone angle less than 2}, J. Amer. Math. Soc. {\bf 28} (2015), no.
1, 199Ð234.

\noindent
[CDS3]    ------------- ,  {\sl K¬ahler-Einstein metrics on Fano manifolds. III:
Limits as cone angle approaches  2 and  completion of the main proof},  J. Amer. Math. Soc. {\bf 28} (2015), no.
1, 235Ð278.

\noindent
{[CP]}  M. Cirant and K. Payne, {\sl Comparison principles for viscosity solutions of 
elliptic branches of fully nonlinear equations independent of the gradient}, ArXiv:2001.09658.

\noindent
[CJY] T. Collins, A. Jacob and S.-T.  Yau,  {\sl  (1, 1) forms with specified Lagrangian
phase: A priori estimates and algebraic obstructions,}, arXiv:1508.01934.

\noindent
[CPW] T. Collins, S. Picard and X. Wu, {\sl Concavity of the Lagrangian phase operator and applications},
Calc. Var. and Partial Differential Equations {\bf 56} (2017), no. 4, Art. 89.
 ArXiv:1607.07194.

\noindent
[CXY]  T.C. Collins, D. Xie, and S.-T. Yau,  {\sl  The deformed Hermitian-Yang-Mills equation
in geometry and physics}, arXiv:1712.00893.

\noindent
[CSh]  T.C. Collins and Y. Shi,  {\sl  Stability and the deformed Hermitian-Yang-Mills equation}, arXiv:2004.04831.

  \noindent
  [CS]  A. Clarke and G. Smith, {\sl 
  The Perron Method and the Non-Linear Plateau problem},  Geom. Dedicata, {\bf 163}, no. 1, (2013), 159-165.

 \noindent
[DDT]  S. Dinew, H-S. Do and T. D. T\^o, {\sl A viscosity approach to the Dirichlet problem for
degenerate complex Hessian type equations}, ArXiv:1712.08572.

\noindent
{[De]  M.  Dellatorre, {\sl 
The degenerate special Lagrangian equation on Riemannian manifolds}, International Math. Res. Notices
(to appear).

\noindent
[D] S. K. Donaldson,  {\sl  Moment Maps and diffeomorphisms},
 Sir Michael Atiyah: a great mathematician of the twentieth century. Asian J. Math. 3 (1999), no. 1, 1Ð15.

\noindent
{[F]}   L. Fu,  {\sl
An analogue of BernsteinÕs theorem}, Houston J. Math. {\bf 24} (1998), 415-419.

\noindent
{[G]}   L.  G\aa rding,  {\sl
An inequality for hyperbolic polynomials}, J. Math. Mech. {\bf 8} no. 2 (1959), 957-965.

\noindent
[Gr] M. Gross, {\sl Special Lagrangian fibrations II: Geometry. A survey of techniques in the
study of special Lagrangian fibrations}, Surv. Differ. Geom., 5, Int. Press, Boston, MA,
1999.

 \noindent
{[CG]} F. Reese Harvey and H. Blaine Lawson {\sl  Calibrated geometries}, Acta Mathematica 
{\bf 148} (1982), 47-157.

 \noindent
{[DD]}   ------------- ,   {\sl  Dirichlet duality and the non-linear Dirichlet problem},
   Comm. on Pure and Applied Math. {\bf 62} (2009), 396-443.

  \noindent
 {[DDR]}  ------------- ,  {\sl   Dirichlet duality and the non-linear Dirichlet problem on Riemannian manifolds}, 
   J. Diff. Geom.  {\bf 88} No. 3 (2011), 395-482.  ArXiv:0907.1981.

  \noindent
 {[HP]}  ------------- ,  {\sl  Hyperbolic polynomials and the Dirichlet problem}, 
    ArXiv:0912.5220.

\noindent
[HP2]  ------------- ,    {\sl  G\aa rding's theory of hyperbolic polynomials},
   {\sl Communications in Pure and Applied Mathematics}  {\bf 66} no. 7 (2013), 1102-1128.

\noindent
[Rest]   ------------- ,   {\sl  The restriction theorem for fully nonlinear subequations},
     Ann. Inst.  Fourier, {\bf 64}  no. 1 (2014), p. 217-265.  ArXiv:1101.4850.

\noindent
[LMA]   ------------- ,   {\sl  Lagrangian potential theory and a Lagrangian equation of  Monge-Amp\`ere type},  
 pp. 217- 257 in Surveys in Differential Geometry, 
 edited by H.-D. Cao, J. Li, R.  Schoen and S.-T. Yau, {\bf  22},  International Press, Somerville, MA, 2018.
  ArXiv:1712.03525.

  \noindent
{[IDP]}  ------------- ,  {\sl The inhomogeneous Dirichlet problem for natural operators on manifolds},
    ArXiv:1805.11121.

\noindent
[H]  N. J.  Hitchin,  {\sl
The moduli space of special Lagrangian submanifolds}, 
Ann. Scuola Norm. Sup. Pisa Cl. Sci. (4)  {\bf 25} (1997), no. 3-4, 503Ð515 (1998).

    \noindent
  [J]  A. Jacob,  {\sl  Weak geodesics for the deformed Hermitian-Yang-Mills equation}, arXiv: 1906.07128.
    
  \noindent
[JY] A. Jacob  and S.-T.  Yau,  {\sl   A special Lagrangian type equation for holomorphic
line bundles}, Math. Ann. {\bf 369} (2017), no.1-2, 869-898.

  \noindent
[JX]   J. Jost and  Y.-L. Xin,   {\sl    A Bernstein theorem for special Lagrangian graphs}, Calc. Var. Partial Differential
Equations. {\bf 15} (2002), 299Ð312.

\noindent
[Jo1] D. Joyce, {\sl Conjectures on Bridgeland stability for Fukaya categories of Calabi-Yau manifolds,
special Lagrangians, and Lagrangian mean curvature flow}, EMS Surv. Math.
Sci. 2 (2015), no. 1, 1Ð62.

\noindent
[Jo2]   ------------- ,  {\sl Special Lagrangian 3-folds and integrable systems. Surveys on geometry and integrable systems},  189Ð233, Adv. Stud. Pure Math., 51, Math. Soc. Japan, Tokyo, 2008.

\noindent
[Jo3]   ------------- ,  {\sl Lectures on special Lagrangian geometry}. Global theory of minimal surfaces, 667Ð695, Clay Math. Proc., 2, Amer. Math. Soc., Providence, RI, 2005.

\noindent
[JLS]  D. Joyce, Y.-I. Lee and R. Schoen,  {\sl  On the existence of Hamiltonian stationary Lagrangian submanifolds in symplectic manifolds},  Amer. J. Math. {\bf 133} (2011), no. 4, 1067Ð1092.

\noindent
[LYZ] N. C. Leung,  S.-T. Yau  and  E. Zaslow,  {\sl From special Lagrangian to
Hermitian-Yang-Mills via Fourier-Mukai transform}, Adv. Theor. Math. Phys.
4 (2000), no. 6, 1319-1341.

\noindent
[NV]  N. Nadirashvili,  and  S. Vlùadutü,  {\sl  Singular solution to the Special Lagrangian
Equations}, Ann. Inst. H. Poincar«e Anal. Non Lin«eaire {\bf 27} (2010), no. 5, 1179-1188.

\noindent
[N1]   A. Neves,   {\sl   Singularities of Lagrangian mean curvature flow: zero-Maslov class case}, 
 Invent. Math. {\bf 168} (2007), no. 3, 449Ð484.

\noindent
[N2]      ------------- ,     {\sl Recent progress on singularities of Lagrangian mean curvature flow}, 
 Surveys in geometric analysis and relativity, 413Ð438, Adv. Lect. Math. (ALM), 20, Int. Press, Somerville, MA, 2011.

\noindent
[N3]      ------------- ,   {\sl   Finite time singularities for Lagrangian mean curvature flow},  
Ann. of Math. (2) {\bf 177} (2013), no. 3, 1029Ð1076.

\noindent
[RS]  Y. Rubinstein and J. Solomon, {\sl The degenerate special Lagrangian equation,}
Adv. Math. 310 (2017), 889-939.

\noindent
[GS1]  G. Smith, {\sl
Special Lagrangian curvature},   
Math. Ann. {\bf 355} (2013), no. 1, 57Ð95.

\noindent
[GS2]   ------------- ,  {\sl The non-linear Dirichlet problem in Hadamard manifolds}, ArXiv:0908.3590.

\noindent
[GS3]   ------------- ,  {\sl The non-linear Plateau problem in non-positively curved manifolds},
 Trans. Amer. Math. Soc., {\bf 365}, (2013), 1109-1124.

\noindent
[GS4]   ------------- ,  {\sl The Plateau problem for convex curvature functions,}, 
to appear in Ann. Inst. Fourier.   ArXiv:1008.3545

\noindent
[SW1]  R. Schoen and J. Wolfson, {\sl Minimizing area among Lagrangian surfaces: the mapping
problem},  J. Differential Geom. {\bf 58} (2001), 1-86.

\noindent
[SW2]   ------------- ,  {\sl
Minimizing volume among Lagrangian submanifolds}. Differential equations: La Pietra 1996 (Florence), 181Ð199, Proc. Sympos. Pure Math., 65, Amer. Math. Soc., Providence, RI, 1999.

\noindent
[SYZ] A.  Strominger,  S.-T.  Yau and  E.  Zaslow,  {\sl Mirror symmetry is T -duality},
Nuclear Phys. B 479 (1996), no. 1-2, 243-259.

\noindent
[T] R. Takahashi,  {\sl Tan-concavity property for Lagrangian phase operators and applications to 
the tangent Lagrangian phase flow},
ArXiv:2002.05132.

\noindent
[W]  M.-T.  Wang,   {\sl Some recent developments in Lagrangian mean curvature flows,
Surveys in Differential Geometry}, Vol. X II. Geometric flows, Int. Press,
Somerville, MA, 2008.

\noindent
[WY1] D. Wang and Y.  Yuan,  {\sl Singular solutions to the special Lagrangian equations
with subcritical phases and minimal surface systems}, Amer. J. Math {\bf 135} (2013),
no. 5, 1157-1177.

\noindent
[WY2]  ------------- ,  {\sl Hessian estimates for special Lagrangian equations
with critical and supercritical phases in general dimensions}, Amer. J. Math.,
{\bf 136} (2014), 481-499.

\noindent
[Wa]  M. Warren,  {\sl Calibrations associated to Monge-Amp`ere equations}, Trans. Amer. Math. Soc.
{\bf 362} (2010), 3947Ð3962.

\noindent
[WaY] Warren, M., and Yuan, Y. Hessian and gradient estimates for three dimensional
special Lagrangian equations with large phase, Amer. J. Math., 132 (2010), 751-
770.

\noindent
[Wo]  J.Wolfson,  {\sl Lagrangian homology classes without regular minimizers},  J. Differential
Geom.  {\bf 71}  (2005), 307-313.

\noindent
[Y0]  
 Yuan, Yu, {\sl A priori estimates for solutions of fully nonlinear special Lagrangian equations},  Ann.\ Inst.\ H.\ Poincar\'e Anal.\ Non Lin\'eaire {\bf 18} (2001), no.\ 2, 261-270.

\noindent
[Y1]   Yu Yuan, {\sl A Bernstein problem for special Lagrangian equations}, Invent. Math. {\bf 150} (2002), no. 1, 117Ð125

  \noindent
{[Y2]}  ------------- , {\sl Global solutions to special Lagrangian equations},  Proc.  A.M.S. {\bf 134} no. 5 (2005),
1355-1358.

\smallskip

\end{document}